\documentclass[11pt]{article}
\usepackage{amsmath,amssymb,mathrsfs,amsfonts}
\usepackage{hyperref}
\usepackage{sectsty}
\input xy
\xyoption{all}

\subsectionfont{\mdseries\itshape}

\subsubsectionfont{\mdseries\itshape}
\newtheorem{Thm}{Theorem}[section]
\newtheorem{Prop}[Thm]{Proposition}
\newtheorem{Lemma}[Thm]{Lemma}
\newtheorem{Cor}[Thm]{Corollary}

\newcommand{\pf}{\noindent{\em Proof.}\ }
\newcommand{\qed}{\hfill $\Box$}

\newcommand{\noin}{\noindent}

\newcommand{\ZZ}{\mathbb{Z}}
\newcommand{\RR}{\mathbb{R}}
\newcommand{\CC}{\mathbb{C}}
\newcommand{\PP}{\mathbb{P}}
\renewcommand{\AA}{\mathbb{A}}
\newcommand{\HH}{\mathbb{H}}
\newcommand{\tensor}[1]{\underset{#1}{\otimes}}

\newcommand{\isoto}{\xrightarrow{\sim}}

\newcommand{\Gr}{\operatorname{Gr}}
\newcommand{\Homsh}{\underline{\mathrm{Hom}}}
\newcommand{\Image}{\operatorname{Image}}
\newcommand{\tot}{\operatorname{tot}}

\newcommand{\OO}{\mathcal{O}}
\newcommand{\RD}{\mathcal{R}}

\newcommand{\om}{\omega}
\newcommand{\Om}{\Omega}
\newcommand{\OL}{\breve{\Om}}
\newcommand{\nb}{\nabla}
\newcommand{\pt}{\partial}
\newcommand{\eps}{\varepsilon}
\newcommand{\lam}{\lambda}
\newcommand{\lag}{\langle}
\newcommand{\rag}{\rangle}
\newcommand{\gen}[1]{\lag #1 \rag}
\newcommand{\ggen}[1]{\lag\hspace{-2pt}\lag #1 \rag\hspace{-2pt}\rag}
\newcommand{\flr}[1]{\lfloor #1 \rfloor}
\newcommand{\ceil}[1]{\lceil #1 \rceil}
\newcommand{\wt}{\widetilde}
\newcommand{\red}[1]{(#1)_{\mathrm{red}}}
\newcommand{\posz}[1]{\underset{\blacktriangle}{{#1}}}
\newcommand{\Disc}{\mathbb{D}}
\newcommand{\HdR}{H_{\mathrm{dR}}}
\newcommand{\HdRc}{H_{\mathrm{dR,c}}}
\newcommand{\cpt}[1]{#1_{\mathrm{c}}}

\begin{document}

\title{Irregular Hodge Filtration
on Twisted De Rham Cohomology\footnote{
This work was partially supported by
the Golden-Jade fellowship of Kenda Foundation,
the NCTS (Main Office) and the NSC, Taiwan.}}
\author{\sc Jeng-Daw Yu\footnote{
Department of Mathematics,
National Taiwan University,
Taipei, Taiwan;
{\tt jdyu@math.ntu.edu.tw}.}}
\maketitle

\begin{abstract}
We give a definition and study the basic properties
of the irregular Hodge filtration
on the exponentially twisted de Rham cohomology
of a smooth quasi-projective complex variety.
\end{abstract}
\section*{Introduction}

\subsection{The goal}
Let $U$ be a complex smooth quasi-projective variety
and $f \in \OO(U)$ be a regular function on $U$.
Consider the algebraic connection $\nb = \nb_f = d+df$
on the structure sheaf $\OO$ of $U$ defined by
\begin{eqnarray*}
\nb : \OO &\to& \Om^1 \\
	v &\mapsto& dv + v\cdot df.
\end{eqnarray*}
It is clear that $\nb$ is integrable
and hence extends to a chain map, still denoted by $\nb$,
on the sheaves $\Om^\bullet$ of differential forms of $U$.
The hypercohomology of the complex $(\Om^\bullet,\nb)$ on $U$
is by definition the de Rham cohomology $\HdR(U,\nb)$
of the connection $\nb$,
which is a finite collection of finite dimensional complex vector spaces.
When $f$ is a constant,
we recover the algebraic de Rham cohomology $\HdR(U/\CC)$ of $U$,
which is equipped with a Hodge filtration
coming from various truncations $(\Om^{\geq p},d)$
of the usual de Rham complex.

($\nb$ is the \textit{exponential twist} of the usual differential $d$
in the sense that the diagram
\[\xymatrix{ \OO \ar[r]^\nb\ar[d]_\cong & \Om^1 \ar[d]^\cong \\
\exp f \cdot \OO \ar[r]^d & \exp f \cdot \Om^1 }\]
commutes.
Here the vertical arrows are the multiplication by
the exponential $\exp(f)$ of $f$.
However since the function $\exp(f)$ is transcendental if $f$ is non-trivial,
one should regard the $\exp(f)$ in the lower corners as a symbol
and it behaves as the exponential
when taking the differentiation.)
\medskip

When $U$ is a curve,
Deligne \cite[pp.109-128]{DMR},
motivated by the analogues
between algebraic connections with irregular singularities
and lisse \'etale sheaves with wild ramifications,
has defined an \textit{irregular Hodge filtration}
$F^\lam$, indexed by $\lam\in\RR$,
on $\HdR(U,\nb)$.
More precisely,
let $X$ be the smooth compactifiaction of $U$
and $S = X\setminus U$ the complement.
The function $f$ on $U$
then extends to a rational function on $X$
and hence $\nb$ defines a meromorphic connection
\[ \nb: \OO_X(*S) \to \Om_X^1(*S) \]
between functions and forms with poles supported on $S$.
We have
\begin{equation}\label{Eq:Intro-U-X}
\HdR^i(U,\nb) =
	\HH^i\left(X,\OO_X(*S) \xrightarrow{\nb} \Om_X^1(*S)\right).
\end{equation}
Deligne then defines an exhaustive and separated decreasing filtration
$F^\lam(\nb)$ on the above two-term complex.
The desired irregular Hodge filtration on $\HdR$
is then given by
\begin{equation}\label{Eq:Intro-Image}
F^\lam\HdR^i(U,\nb) := \Image\left\{
	\HH^i\left(X,F^\lam(\nb)\right) \to \HdR^i(U,\nb) \right\}
\end{equation}
under the identification \eqref{Eq:Intro-U-X}.
When $f$ is a constant,
$F^\bullet(\nb)$ reduces to the pole-order filtration $P^\bullet$
defined in \cite[II.3.12]{Del163}
and thus one recovers the usual Hodge filtration.

Deligne has shown that
the spectral sequence associated with this filtration
degenerates at the initial stage
(i.e.~the arrow in \eqref{Eq:Intro-Image} is always injective)
and proved that the irregular Hodge filtration
respects the pairing
between the cohomology of $\nb$ and of the dual connection.
However the filtration and its complex conjugate
(with respect to the real structure from its Betti counterpart)
are not opposite to each other in general.
We remark that in \cite{DMR},
the definition of the irregular Hodge filtration is justified
by its relation with
the expected weights of special values of the gamma function.
Moreover the filtration is defined for more general connections
of certain type,
not necessarily of rank one.

\medskip
In this paper,
we propose a definition of the irregular Hodge filtration
for $\nb$ on $U$ of arbitrary dimension.
The idea is similar to Deligne's approach.
We first pick a compactification $X$ of $U$
such that $f$ extends to a morphism from $X$ to $\PP^1$
and that the complement $S:=X\setminus U$
is a normal crossing divisor.
We then define a decreasing filtration $F^\lam(\nb)$
on the twisted meromorphic de Rham complex
$(\Om_X^\bullet(*S),\nb)$.
The sheaves involved in each $F^\lam(\nb)$ are all locally free on $X$.
However the new definition does not coincide with Deligne's
when $U$ is a curve.
In fact when $f$ is constant,
our filtration is not the same as the pole-order filtration $P^\bullet$
but is equal to the usual Hodge filtration
$(\Om_X^{\geq \lam}(\log S),d)$
of the de Rham complex of logarithmic differential forms.
(Thus one still recovers the usual Hodge filtration for $U$
in this case.)
Although the filtration fails to be exhaustive in general,
it is rich enough to capture the de Rham cohomology of $\nb$
and indeed induces the same filtration on the cohomology
as Deligne's in the curve case.

Another advantage of using the logarithmic differential forms
lies in the fact that
unlike the curve case,
one does not have a canonical choice of the compact $X$.
Two different choices are connected by a birational morphism $\pi$
and the sheaves $\Om_X^p(\log S)$ behave well
under $\pi$.
We shall prove that the irregular Hodge filtration
on the de Rham cohomology of $\nb$ obtained in this way
is independent of the choice of $X$
and satisfies some functorial properties.
As in the curve case,
we also demonstrate that the filtration
respects the Poincar\'e pairing
between the cohomology of $\nb$ and of the dual connection.

After the appearance in arXiv of the first version of this paper,
we obtain a proof of the $E_1$-degeneracy
of the spectral sequence associated with
the irregular Hodge filtration
in the joint work \cite{ESY} with H.~Esnault and C.~Sabbah.

\medskip
Finally we mention that in the direction
of relating the algebraic de Rham cohomology
to a Betti type cohomology attached to
any integrable algebraic connection
of arbitrary rank
via periods,
the homology with coefficients in
rapid decay simplicial chains has been defined
and the duality to the de Rham cohomology has been established
in \cite{BE} for the curve case
and \cite{Hien} in general.
On the other hand,
the relation to the nonabelian Hodge theory
has been discussed in \cite{OV}.
The exponentially twisted de Rham cohomology
also appears in the theory of
mirror symmetry \cite{KKP}
and the study of Donaldson-Thomas invariants \cite{KS}.
We hope the investigation of the irregular Hodge filtration
can provide more structures and shed some light
into these areas.
In \cite{S_II},
another generalization of the irregular Hodge filtration
in the higher rank case over a projective line
is developed and has been connected to
the so-called supersymmatric index.

\subsection{The structure of the paper}

After the introductory section
we give the definition of the irregular Hodge filtration
of the twisted de Rham complex
on a certain compactification $X$ of $U$
in \S\ref{Sect:IrrHodge}.
We show that the induced filtration
on the de Rham cohomology of $(U,\nb)$
is independent of $X$.
Along the way some basic properties of the filtration are derived.
We shall define the corresponding filtration
on the cohomology with compact support.
We then establish the perfect Poincar\'e pairing between
the de Rham cohomology of $\nb$
and of its dual with compact support
in \S\ref{Sect:Duality}.
The irregular Hodge filtrations on them are shown
to respect the Poincar\'e pairing.
In fact
since we do not know
when the Hodge to de Rham spectral sequence degenerates here,
we will also define pairings between terms on each stage
of the spectral sequence
and discuss their relations.

\S\ref{Sect:Hypersurf} and \S\ref{Sect:Toric}
are devoted to providing examples.
In \S\ref{Sect:Hypersurf}
we discuss the case where $U = \AA^1\times U'$
and $f$ is the direct product of the identity on $\AA^1$
and a function $f'$ on $U'$.
In this case,
the irregular Hodge filtration reduces
to the usual Hodge filtration of the subvariety defined by $f'$
if it is smooth.
We then recall the work of Adolphson and Sperber
on the twisted de Rham cohomology
over a torus in \S\ref{Sect:Toric}.
In this case
a filtration coming from the Newton polyhedron $\Delta$ of the function $f$
is defined
and the associated spectral sequence is shown
to degenerate if $f$ is \textit{non-degenerate} with respect to $\Delta$
in \cite{AS-Nagoya}.
We show that in this case
the filtration from $\Delta$ on the de Rham cohomology
coincides with our irregular Hodge filtration.

Finally in the appendix
we briefly recall Deligne's definition of the irregular Hodge filtration
in the curve case
and indicate that
his definition gives the same filtration
on the de Rham cohomology as ours.

\medskip
Part of this work has been done
during my visit of Universit\"at Duisburg - Essen
in 2011 - 2012.
I am grateful for the financial support,
the hospitality and the inspiring environment
provided by the group of
\textit{Essener Seminar f\"ur Algebraische Geometrie und Arithmetik}.
I thank Professor H\'el\`ene Esnault
for helpful discussions and bringing my attention to the paper \cite{S_II}.
Thanks are also due to the referee
for the careful reading and suggestions
which help improve the presentation.

\subsection{Notations and conventions}
To shorten the notation,
let
\[ \AA = \AA^1 \quad\text{and}\quad \PP=\PP^1 \]
be the affine line and the projective line, respectively
in the rest of this paper.
Let
\[ \Disc \quad\text{and}\quad \Disc^\circ = \Disc\setminus\{0\} \]
be the open unit disc and the punctured disc
of the complex plane, respectively.
For a divisor $D$ on a variety,
$\red{D}$ denotes the associated reduced subvariety,
i.e., the support of $D$.

Since we will use the sheaves of logarithmic differentials intensively,
we introduce the following notation.
Let $X$ be a smooth completion of $U$
such that the complement $S = X\setminus U$ is a normal crossing divisor.
We let
\[ \OL^p = \OL_X^p = \OL^p_{U\subset X} := \Om^p_X(\log S) \]
be the sheaf on $X$ of differential forms of degree $p$,
regular on $U$ and with at worst logarithmic poles along $S$.

For a decreasing filtration $F^\lam$ indexed by $\lam \in \RR$,
we set
\begin{eqnarray*}
F^{\lam-} = \bigcap_{i<\lam} F^i &\text{and}&
	F^{\lam+} = \bigcup_{i>\lam} F^i.
\end{eqnarray*}
The $\lam$-th graded piece $\Gr^\lam = \Gr^\lam_F$ of $F^\bullet$
is defined as $F^\lam/F^{\lam+}$.

For a complex $K = (K^\bullet,\delta^\bullet)$,
the degree $p$ term of the shift $K[n]$
is $K^{n+p}$ with the differential $\delta^{n+p}$.
If we want to locate the degree 0 term of a complex
to avoid confusion,
we put the symbol $\blacktriangle$ under that term,
e.g.,
$\cdots\to\posz{A}\to B\to\cdots$.
The use of $\blacktriangle$ in some variants
in the paper should be clear.
For a double complex $(K^{\bullet,\bullet},\delta_1,\delta_2)$,
the symbol $\tot(K^{\bullet,\bullet})$
denotes the total complex attached to $K^{\bullet,\bullet}$
with differential $\delta_1+(-1)^p\delta_2$ on $K^{p,q}$.

\section{The irregular Hodge filtration}\label{Sect:IrrHodge}

\subsection{The de Rham cohomology
	and good compactifications}
Fix a complex smooth quasi-projective variety $U$
and a global regular function $f$ on it,
regarded as an element $f\in \OO(U)$
or a morphism $f:U\to\AA$ interchangeably.
As in the introduction,
let $\nb = \nb_f = d+df$ be the integrable connection
on the structure sheaf $\OO$ of $U$.
It then extends to the twisted de Rham complex on $U$
\[ (\Om^\bullet,\nb) =
	\left[\OO \xrightarrow{\nb} \Om^1 \xrightarrow{\nb}
		\Om^2 \to\cdots\right]. \]

\bigskip
\noin\textit{Definition.}
The \textit{de Rham cohomology of the connection $\nb$}
is the hypercohomology
\[ \HdR^i(U,\nb) := \HH^i\big(U,(\Om^\bullet,\nb)\big). \]

\bigskip
\noin\textit{Definition.}
Let $j:U\to X$ be a compactification of $U$
with the complement $S:=X\setminus U$.
The pair $(X,S)$ is called
a \textit{good compactification of $(U,f)$}
if $S$ is a normal crossing divisor of $X$
and $f$ extends to a morphism $f:X\to\PP$.
In this case we have the commutative diagram
\[\xymatrix{
U \ar[r]^f\ar[d]_j & \AA \ar[d] \\
X \ar[r]^f & \PP. }\] 
\bigskip

By the elimination of indeterminacy
and the resolution of singularities
there always exists a good compactification $(X,S)$ of $(U,f)$.
Given such an $X$
and a point $a\in f^{-1}(\infty)$ of $X$,
there exists a system of analytically local coordinates
\[ \{x_1,\cdots,x_l,t_1,\cdots,t_m,y_1,\cdots,y_r\}
\quad\text{for some $l,m,r\geq 0$} \]
such that
\begin{itemize}
\item $S = (xt)$ is a union of coordinate hyperplanes, and
\item $f = \frac{1}{x^e}f_0$
for some exponent $e\in\ZZ_{> 0}^l$
and some analytic $f_0$ with $f_0(a)\neq 0$.
\end{itemize}
This local picture will be used repeatedly.
\medskip

On the other hand,
the connection $\nb$ on $U$ extends to the twisted complex
\begin{equation}\label{Eq:dRCpx}
(\Om_X^\bullet(*S),\nb) =
	\left[\OO_X(*S) \xrightarrow{\nb} \Om_X^1(*S) \xrightarrow{\nb}
		\Om_X^2(*S) \to\cdots\right].
\end{equation}
Since $\Om_X^p(*S) = j_*\Om_U^p \isoto \RD j_*\Om_U^p$,
we have
\begin{equation}\label{Eq:dR-U-X}
\HdR^i(U,\nb) = \HH^i\big(X,(\Om_X^\bullet(*S),\nb)\big).
\end{equation}

\subsection{The Hodge filtration on the de Rham complex}
Fix a good compactification $(X,S)$ of $(U,f)$.
We shall define on the complex \eqref{Eq:dRCpx}
a separated filtration $F^\lam$,
indexed by $\lam \in \RR$,
which is left continuous
(i.e., $F^\lam = F^{\lam-}$)
and with discrete jumps
(i.e.~the set $\{\lam\in\RR \,|\, \Gr^\lam \neq 0\}$
is discrete).
(It will also be exhaustive if $f:U\to\AA$ is proper.)

Let $P$ be the pole divisor of $f$ on $X$;
it is effective and supported on $S$.
We have
\[ f\in\OO_X(P) \quad\text{and}\quad df\in\OL_X^1(P). \]
(Recall that $\OL_X^p := \Om_X^p(\log S)$.)

\bigskip
\noin\textit{Definition.}
Let
\begin{equation}\label{Eq:F0lam}
F^0(\lam) :=
\left[ \OO(\flr{-\lam P}) \xrightarrow{\nb} \OL^1(\flr{(1-\lam) P})
	\to\cdots\to \OL^p(\flr{(p-\lam) P}) \to\cdots \right],
\end{equation}
regarded as a subcomplex of \eqref{Eq:dRCpx}.
The \textit{irregular Hodge filtration of $\nb$}
is the filtration on \eqref{Eq:dRCpx} defined by
\[ F^\lam(\nb) = F^0(\lam)^{\geq\ceil{\lam}}
	\quad(\lam\in\RR). \]
We use $F^\lam(\nb)^p$ to denote the degree $p$ component
of $F^\lam(\nb)$;
it is a locally free subsheaf of $\Om_X^p(*S)$.
\bigskip

Clearly at degree $p$, we have
\begin{equation}\label{Eq:Ind-Fil-I}
F^\lam(\nb)^p = \left\{\begin{array}{ll}
	0 & \text{if $p < \lam$} \\
	\OL^p(\flr{(p-\lam)P}) & \text{if $p \geq \lam$} \end{array}\right.
\end{equation}
and that $F^\bullet(\nb)$ obeys the following two rules:
\begin{equation}\label{Eq:Ind-Fil-II}
\begin{array}{rcll}
F^\lam(\nb)^0 &=& \left(F^{\lam+1}(\nb)^0\right)(P)
	& \text{if $\lam \leq -1$}, \\ \\
F^\lam(\nb)^p &=&
	\OL^p \tensor{\OO_X} F^{\lam - p}(\nb)^0
	& \text{for all $\lam \in \RR$}.
\end{array}\end{equation}

In the rest of this subsection,
we build up some basic properties of this filtration.
\medskip

First consider the local situation.
Let $U = (\Disc^\circ)^l\times(\Disc^\circ)^m\times\Disc^r$
with coordinates
\[ \{x_1,\cdots,x_l,t_1,\cdots,t_m,y_1,\cdots,y_r\}. \]
Let $f = \frac{f_0}{x_1^{e_1}\cdots x_l^{e_l}}$
with $e_i>0$,
$f_0$ regular and nowhere vanishing on $\Disc^{l+m+r}$,
and $\nb = \nb_f$ the associated connection on $U$.
Let $\wt{U} = U\times\Disc^a\times(\Disc^\circ)^b$
with the natural embedding into $X = \Disc^{l+m+r+a+b}$
and $S:=X\setminus\wt{U}$.
On $X$,
let $F^\lam(\wt{\nb})$ be the filtration of the connection $\wt{\nb}$
attached to $f$ regarded as a function on $\wt{U}$,
and
$F_\boxtimes^\lam$ the exterior product filtration
of $F^\mu(\nb)$ and $F^\nu(d)$.
One checks directly that
$F_\boxtimes^\lam$ is a subcomplex of $F^\lam(\wt{\nb})$.

\begin{Prop}\label{Prop:Prod-Fil}
In the local setting as above, the natural inclusion
\[ F_\boxtimes^\lam \,\,\to\, F^\lam(\wt{\nb}) \]
of subcomplexes of $(\Om_X^\bullet(*S),\nb)$
is a quasi-isomorphism for each $\lam$.
\end{Prop}

\pf
We first consider the case where $a=1$ and $b=0$.
Let $n = l+m+r$.
Fix $\lam\leq n+1$.
The quotient of the natural inclusion of complexes
\[\xymatrix{
F^\lam_\boxtimes \left( \big[\OO_{\Disc^n}(*S) \xrightarrow{\nb} \Om^1(*S)
	\to\cdots\to \Om^n(*S)\big]
	\boxtimes \big[\OO_\Disc \xrightarrow{d} \Om^1\big] \right) \ar[d] \\
F^\lam \left(\OO_{X}(*S) \xrightarrow{\nb} \Om^1(*S)
	\to\cdots \to \Om^{n+1}(*S)\right) }\]
is described as follows.
Let $z$ be the coordinate of the last piece $\Disc$ of $\wt{U}$.
Let $A = \OO(\wt{U})$
and $\Psi =$ all possible exterior products among the 1-forms in
$\left\{\frac{dt_i}{t_i},dy_j\right\}$ of degree $\geq \lam$.
Then, as an $A$-module, the quotient decomposes into
\[ \bigoplus_{\om\in\Psi} \left( B(\om),\nb|_{B(\om)} \right) \]
with $B(\om)[p+1] =$
\[\begin{split}
\left(\frac{1}{x^{\flr{(p+1-\lam)e}}}A\Big/\frac{1}{x^{\flr{(p-\lam)e}}}A\right)
	&\om_0 dz
	\xrightarrow{\text{left multip.~by $dx^{-e}$}} \\
	\bigoplus_{i=1}^r &
\left(\frac{1}{x^{\flr{(p+2-\lam)e}}}A\Big/\frac{1}{x^{\flr{(p+1-\lam)e}}}A\right)
	\frac{dx_i}{x_i}\om_1 dz \\
	\to &\bigoplus_{1\leq i<j\leq r}
\left(\frac{1}{x^{\flr{(p+3-\lam)e}}}A\Big/\frac{1}{x^{\flr{(p+2-\lam)e}}}A\right)
	\frac{dx_i}{x_i}\frac{dx_j}{x_j}\om_2 dz
	\to\cdots \end{split}\]
where $p = \deg(\om)$
and $\om_k = f_0^k\cdot\om$.
It is clear that
the complex $B(\om)[p+1]$ of $\CC$-vector spaces is isomorphic to
\[ \left(\frac{1}{x^{\flr{(p+1-\lam)e}}}A\Big/\frac{1}{x^{\flr{(p-\lam)e}}}A\right)
	\tensor{\CC}
	\left[\bigwedge^0 C
	\xrightarrow{\text{left multip.~by $\sum_{i=1}^l v_i$}}
		\bigwedge^1 C \to\cdots\to \bigwedge^l C \right] \]
where the later is the total Koszul complex attached to
the $\CC$-vector space $C$ generated by the basis
$\left\{v_i = -e_i\frac{dx_i}{x_i}\right\}_{i=1}^l$.
Now this Koszul complex has null-cohomology
and thus the assertion follows in this case.

For the case where $a=0$ and $b=1$,
one simply replaces $dz$ by $\frac{dz}{z}$ in the above arguments.

The general case then follows from the above two cases inductively
and the fact that
the usual Hodge filtration of the logarithmic de Rham complex of
$\Disc^a\times(\Disc^\circ)^b$
is equal to the product filtration of the filtrations on its factors.
\qed

\begin{Prop}\label{Prop:Change-Support}
Let $D$ and $E$ be divisors of $X$ supported on $S$ and $\red{P}$,
respectively.
Suppose $E$ is effective.
Then the natural inclusion
\begin{equation}\label{Eq:add-divisor}
\vcenter{\xymatrix{
\left[\OO(D) \xrightarrow{\nb} \OL^1(D+P) \to\cdots\to \OL^p(D+pP)
	\to\cdots \right]
	\ar[d] \\
\left[ \OO(D+E) \xrightarrow{\nb} \OL^1(D+E+P) \to\cdots\to
	\OL^p(D+E+pP) \to\cdots \right] }}
\end{equation}
of complexes on $X$ is a quasi-isomorphism.
\end{Prop}

Indeed by induction, it suffices to consider the case
where $E$ is an irreducible component of $\red{P}$.
The assertion is then obtained
by a local computation similar to the proof of Prop.\ref{Prop:Prod-Fil}.
We omit the details.

\begin{Cor}\label{Cor:F0-F}
The inclusion
\[ \Big(F^0(\nb)\,\text{with the induced irregular Hodge filtration}\Big) \to
	\left(\big(\Om_X^\bullet(*S),\nb\big), F^\lam(\nb)\right) \]
is a quasi-isomorphism of filtered complexes on $X$.
\end{Cor}

\pf
Write $S=\red{P}+T$.
Prop.\ref{Prop:Prod-Fil} (plus \cite[proof of Prop.II.3.13]{Del163})
and the above proposition
give respectively the two quasi-isomorphisms
\begin{eqnarray*}
F^0(\nb) &\isoto& F^0(\nb)(*T) \\
	&\isoto& \left(F^0(\nb)(*T)\right)(*\red{P}),
\end{eqnarray*}
both compatible with the equipped filtrations;
the last term is simply the complex $(\Om_X^\bullet(*S),\nb)$
with the filtration $F^\lam(\nb)$.
\qed
\bigskip

Notice that the corollary above implies immediately
that $\HdR^i(U,\nb)$ is finite dimensional for any $i$
and is zero unless $0\leq i\leq 2\cdot\dim U$,
since it is the hypercohomology
of a chain of coherent sheaves on a compact $X$ of length $\dim U$.
\bigskip

\noin\textit{Definition.}
On $X$
the \textit{logarithmic complex attached to $\nb$}
is the sub-filtered complex
$(\Om_X^\bullet(\log\nb), F^\lam) \subset (\Om_X^\bullet(*S), F^\lam(\nb))$
defined as
\[ \Om_X^\bullet(\log\nb) = F^0(\nb)\tensor{\OO_X}\OO_X(-\red{P}). \]
\bigskip

Inside this logarithmic complex,
$\Om_X^0(\log\nb) = \OO(-\red{P})$ is pure of filter degree 0
while $\Om_X^{>0}(\log\nb)$ is of positive filter degree,
i.e.~jumps $>0$.
The following corollary
gives us the information of $\Gr^0(\nb)$.

\begin{Cor}
The inclusion
$\left(\Om_X^\bullet(\log\nb),\nb, F^\lam\right) \to
	\left(\Om_X^\bullet(*S),\nb, F^\lam(\nb)\right)$
is a quasi-isomorphism of filtered complexes.
In particular
we have the quasi-isomorphism
\[ \OO(-\red{P}) \isoto \Gr^0(\nb). \]
\end{Cor}

\pf
This follows from Prop.\ref{Prop:Change-Support}
(by taking $D=-\red{P}$, $E=\red{P}$)
together with the above corollary.
\qed

\subsection{The Hodge filtration on the de Rham cohomology}
In the previous subsection,
we defined the irregular Hodge filtration on the twisted de Rham complex
upon a chosen good compactification $(X,S)$ of $(U,f)$.
Here we prove that the induced filtration on $\HdR(U,\nb)$
does not depend on the choice of $X$.
\medskip

We begin by considering a map $\pi: (X',S') \to (X,S)$
between two good compactifications of $(U,f)$.
The corresponding irregular Hodge filtrations on them
will be denoted by $F_X^\lam(\nb)$ and $F_{X'}^\lam(\nb)$,
respectively.
Recall that
since $\pi$ is a proper birational morphism between smooth varieties,
we have
\[ \RD\pi_*\OO_{X'} = \OO_{X}. \]

\begin{Prop}
With notations as above,
we have $\pi^*F_X^\lam(\nb) \subset F_{X'}^\lam(\nb)$
where $\pi^*$ denotes the componentwise
pullback to $\OO_{X'}$-modules.
In particular we obtain
\[ \pi^*: \HH\left(X,F_X^\lam(\nb)\right) \to
	\HH\left(X',F_{X'}^\lam(\nb)\right). \]
\end{Prop}

\pf
Let $P$ and $P'$ be the pole divisors of $f$
on $X$ and $X'$, respectively.
One sees readily that
\[ \pi^*\OO_X(\flr{\eta P}) \subset \OO_{X'}(\flr{\eta P'})
	\quad\text{for any $\eta\geq 0$}. \]
Since $\pi^*\OL_X^p \subset \OL_{X'}^p$,
the assertion follows from
the identities in \eqref{Eq:Ind-Fil-I} and \eqref{Eq:Ind-Fil-II}.
\qed

\begin{Lemma}\label{Lemma:Push-forward}
Let $(X,S)$ be a good compactification of $(U,f)$.
Suppose that $(X',S')$ is another good compactification
obtained by a blowup $\pi: X' \to X$ along a smooth center
which has normal crossing with $S$.
Then we have the following.
\begin{enumerate}
\item The adjunction map $\OL_X^p \to \RD\pi_*\OL_{X'}^p$
of the pullback $\pi^*\OL_X^p \to \OL_{X'}^p$
is a quasi-isomorphism for any $p\in\ZZ$.
\item The adjunction map
\[ F_X^\lam(\nb) \to \RD\pi_*\left(F_{X'}^\lam(\nb)\right) \]
of the natural inclusion $\pi^*F_X^\lam(\nb) \to F_{X'}^\lam(\nb)$
is a quasi-isomorphism for any $\lam\in\RR$.
\end{enumerate}
\end{Lemma}

\pf
(i) At a point in the center $\Xi$ of the blowup,
there exist local coordinates $\{x_1,\cdots, x_n\}$ of $X$
and three positive integers $r, a,b$ with $a\leq r\leq b\leq n$
such that
$S = (x_1\cdots x_r)$
and $\Xi$ is defined by
\[\left\{\begin{array}{ll}
x_i=0 & \text{if $1\leq i\leq a$} \\
x_j=0 & \text{if $r<j\leq b$}. \end{array}\right.\]
Let $E$ be the exceptional divisor.
Over this local chart, we have that
$E/\Xi$ is fibered by projective spaces
of dimension $\dim_\Xi E = (a+b-r-1)$.
Then using the standard affine cover of the blowup,
one checks directly that the quotient
$\OL_{X'}^p/\pi^*\OL_X^p$
is isomorphic to
\begin{equation}\label{Eq:Pullback-OL}
\bigoplus_{i=1}^{b-r} \OO_{E/\Xi}(-i)^{\delta_i}
\quad\text{where
$\delta_i = \binom{b-r}{i}\binom{n-r+b}{p-i}$}.
\end{equation}
Since $\RD\pi_*\OO_{E/\Xi}(-i) = 0$
for $1\leq i\leq b-r$,
the assertion follows from the projection formula
(\cite[Exer.III.9.8.3]{Hart} or \cite[Prop.II.5.6]{H-RD}).

(ii) Fix a non-positive real number $\lam$.
Let $Q_\lam$ be the coherent sheaf on $X'$
defined by the short exact sequence
\begin{equation}\label{eqn:diff-blowup}
0 \to \pi^*\left(F^\lam_X(\nb)^0\right) \to F^\lam_{X'}(\nb)^0 \to Q_\lam \to 0.
\end{equation}
Then $Q_\lam$ is concentrated on $E$.
We use the same local coordinates $\{x_1,\cdots, x_n\}$ of $X$
as in (i).
By shrinking the neighborhood if necessary,
we have $f = (x_1^{e_1}\cdots x_r^{e_r})^{-1}f_0$
with $f_0$ regular and nowhere vanishing.
Over this local chart, we have that
\begin{itemize}
\item $\pi^*f$ has pole order $e:=e_1+\cdots +e_a$ along $E$, and
\item above the origin of $X$,
the sequence \eqref{eqn:diff-blowup} is given by
\begin{equation}\label{Eq:Excep-sheaf}
0 \to \OO\left(\flr{-\lam \wt{P}}+\sum_{i=1}^a\flr{-\lam e_i}E\right)
	\to \OO\left(\flr{-\lam \wt{P}}+\flr{-\lam e}E\right) \to Q_\lam \to 0
\end{equation}
where $\wt{P}\subset X'$ denotes the proper transform of $P\subset X$.
(Thus $\red{P'} = \red{\wt{P}}+E$.)
\end{itemize}

Inserting the intermediate locally free sheaves of $X'$ into the inclusion
\[ \OO\left(\flr{-\lam \wt{P}}+\sum_{i=1}^a\flr{-\lam e_i}E\right)
	\subset \OO\left(\flr{-\lam \wt{P}}+\flr{-\lam e}E\right) \]
by adding one more copy of the divisor $E$ in each step,
we get a filtration in the middle term of \eqref{Eq:Excep-sheaf}.
It then induces a filtration on $Q_\lam$.
To get information of the induced grading on $Q_\lam$,
one has to compute the restriction
\[ \OO_E\left(\flr{-\lam \wt{P}}+\flr{-\lam e}E\right) \]
of the sheaf to $E$.
Write $\red{\wt{P}} = \sum_{i=1}^r \wt{P}_i$
where $\wt{P}_i$ is the proper transform
of the $i$-th coordinate hyperplane.
We notice that, still over the origin of $X$,
\[ E.\wt{P}_i = \left\{\begin{array}{cl}
H & \text{if $1\leq i\leq a$} \\
0 & \text{if $i>a$} \end{array}\right.
	\quad,\quad E.E = -H
	\quad(\text{$H$ denotes a hyperplane section}) \]
and
\[ 0\leq \flr{-\lam e} - \sum_{i=1}^a\flr{-\lam e_i} \leq a-1. \]
Therefore
$Q_\lam$ over the origin of $X$
is a successive extension of various
$\OO_{E/\Xi}(-\mu H)$ with $0 < \mu \leq a-1$.

Now together with \eqref{Eq:Pullback-OL},
one then obtains that over the origin of $X$,
the quotient of $\OL_{X'}^p\otimes F^\lam_{X'}(\nb)^0$
by $\pi^*\left(\OL_X^p\otimes F^\lam_X(\nb)^0\right)$
is a successive extension of various
$\OO_{E/\Xi}(-\mu H)^{m_\mu}$
for some $m_\mu\geq 0$
with $0 < \mu \leq (a-1)+(b-r) = \dim_\Xi E$.
With $\mu$ in this range, we have that
$\RD\pi_*\OO_{E/\Xi}(-\mu H)$ is quasi-isomorphic to zero.
By the second identity in \eqref{Eq:Ind-Fil-II}
and the projection formula,
we obtain the stated result.
\qed

\begin{Thm}
The hypercohomology $\HH\left(X, F^\lam(\nb)\right)$
only depends on $(U,f)$,
not on the choice of the good compactification $(X,S)$.
\end{Thm}

\pf
First suppose that $\pi:X'\to X$
is a morphism between good compactifications.
Recall the weak factorization theorem of birational morphisms
\cite[Thm.0.0.1]{W}:
The birational morphism $\pi$ admits a factorization
into the following commutative diagram of birational morphisms
\[\xymatrix{
X' = X_0\ar@{--}[r]\ar@/_2pc/[rrrrrr]_\pi &
	X_1 \ar@{--}[r]\ar@/_1pc/[l] &
	\cdots\ar@{--}[r]_{\alpha_m} &
	X_m\ar@{--}[r]\ar@/_3pc/[lll]_{\beta_m}\ar@/^3pc/[rrr]^{\gamma_m} &
	\cdots\ar@{--}[r]_{\alpha_{n-1}} &
	X_{n-1} \ar@{--}[r]\ar@/^1pc/[r] & X_n = X. }\]
Here, for $1\leq i\leq n$,
\begin{itemize}
\item $X_i$ is a smooth completion of $U$
with $S_i := X_i\setminus U$ a normal crossing divisor;
\item $X_{i-1} \overset{\alpha_i}{---} X_i$ represents
either a blowup $\alpha_i: X_{i-1} \to X_i$
of $X_i$ along a smooth center
which is of normal crossing with $S_i$,
or a blowup $\alpha_i: X_i \to X_{i-1}$
of $X_{i-1}$ along a smooth center
which is of normal crossing with $S_{i-1}$;
\item there exists an integer $m \in [1,n]$ such that
$X_i$ are equipped with morphisms
\[\begin{array}{ll}
\beta_i: X_i \to X_0 = X', & 1\leq i\leq m \\
\gamma_i: X_i \to X_n = X, & m\leq i\leq n. \end{array} \]
\end{itemize}
The first and the third conditions ensure that
each $(X_i,S_i)$ is a good compactification of $(U,f)$.
We let $F_i^\lam(\nb)$ denote the associated
irregular Hodge filtration on $X_i$.

Set $\gamma_0 = \pi$
and $\gamma_i = \pi\circ\beta_i: X_i \to X$ for $1\leq i\leq m$.
Then for each $1\leq i\leq n$,
we have the commutative diagrams
\[\xymatrix{ \text{(I)} & X_{i-1}\ar[rr]^{\alpha_i}\ar[dr]_{\gamma_{i-1}} &&
	X_i\ar[dl]^{\gamma_i} \\
&& X } \qquad\text{or}\qquad
\xymatrix{ \text{(II)} & X_{i-1}\ar[dr]_{\gamma_{i-1}} &&
	X_i\ar[ll]_{\alpha_i}\ar[dl]^{\gamma_i} \\
&& X. }\]
By Lemma \ref{Lemma:Push-forward} (applied to $\alpha_i$),
we obtain
\[ \RD{\gamma_{i-1}}_*F^\lam_{i-1}(\nb) = \left\{\begin{array}{c}
	\RD{\gamma_i}_*\left(\RD{\alpha_i}_*F^\lam_{i-1}(\nb)\right)
		\quad\text{in case (I)} \\ \\
	\RD{\gamma_{i-1}}_*\left(\RD{\alpha_i}_*F^\lam_i(\nb)\right)
		\quad\text{in case (II)} \end{array}\right\}
	= \RD{\gamma_i}_*F^\lam_i(\nb). \]
(The $=$ means quasi-isomorphic.)
Thus by induction on the index $i$ in $\gamma_i$,
one obtains that
$F_X^\lam(\nb) \to \RD\pi_*\left(F_{X'}^\lam(\nb)\right)$
is a quasi-isomrphism.
Therefore the assertion follows in this case.

Now given two good compactifications
$(X_1,S_1)$ and $(X_2,S_2)$ of $(U,f)$,
one can always find a third one that dominates the two.
Indeed we have the standard commutative diagram:
\[\xymatrix{
	& X_1 & \\
U \ar@{^(->}[ur]\ar@{^(->}[dr]\ar@{^(->}[rr] &&
	X \ar[r]\ar[ul]_{\pi_1}\ar[dl]^{\pi_2}
	& \overline{U} \subset X_1 \times X_2 \\
	& X_2 & }\]
where $\overline{U}$ is the closure of $U$ in $X_1\times X_2$
via the diagonal embedding
and $X\to\overline{U}$ is a certain sequence of blowups
such that $(X,X\setminus U)$ is a good compactification.
The above discussion then shows that
$\pi_1$ and $\pi_2$ induce isomorphisms
on the hypercohomology of the corresponding $F^\lam(\nb)$.
This completes the proof.
\qed\bigskip

Applying the snake lemma
to the long exact sequence associated with
\[ 0 \to F^{\lam+}(\nb) \to F^\lam(\nb) \to \Gr^\lam(\nb) \to 0, \]
the above theorem then yields the following.

\begin{Cor}\label{Cor:HGr-well-defined}
The hypercohomology
$\HH^i\left(X, \Gr^\lam(\nb)\right)$
does not depend on the choice of $X$.
\end{Cor}

\noin\textit{Definition.}
Let $(X,S)$ be a good compactification of $(U,f)$.
\begin{enumerate}
\item For any $\lam\in\RR$, we define
\begin{eqnarray*}
H^i\left(U,F^\lam(\nb)\right) &:=& \HH^i\left(X,F^\lam(\nb)\right) \\
H^i\left(U,\Gr^\lam(\nb)\right) &:=& \HH^i\left(X,\Gr^\lam(\nb)\right).
\end{eqnarray*}
\item
The \textit{irregular Hodge filtration} $F^\lam$ on $\HdR^i(U,\nb)$
is defined by setting
\begin{eqnarray*}
F^\lam \HdR^i(U,\nb_f) &=&
	\Image\left\{\HH^i\left(X,F^\lam(\nb)\right) \to
		\HH^i\Big(X,(\Om_X^\bullet(*S),\nb)\Big) \right\} \\
	&=& \Image\left\{H^i\left(U,F^\lam(\nb)\right) \to
		\HdR^i(U,\nb) \right\}
\end{eqnarray*}
induced from the inclusion $F^\lam(\nb) \to (\Om_X^\bullet(*S),\nb)$
and via the canonical isomorphism \eqref{Eq:dR-U-X}.
\end{enumerate}
\bigskip

The definition does not depend on the choice of $X$.
Notice that by Cor.\ref{Cor:F0-F} we have
\[ H^i\left(U,F^\lam(\nb)\right) = \HdR^i(U,\nb)
	\quad\text{if $\lam\leq 0$}. \]

\subsection{The cohomology with compact support}
In the last part of this section
we introduced the de Rham cohomology with compact support
of the connection $\nb$
and define the corresponding irregular Hodge filtration.
For the classical case,
see \cite[\S 4.3]{EZ}.
Again the definitions rely on
choosing a good compactification $X$ first.
It is possible to establish the corresponding properties
for the cohomology with compact support
and prove that the definition of the irregular Hodge filtration
does not depend on the choice of $X$
as in the previous discussion.
However we do not proceed in this direction.
The independency will be clear
once we obtain the duality in the next section.
Notice that the proofs of the results in the next section
do not use the proposition below.

\bigskip
\noin\textit{Definition.}
Let $(X,S)$ be a good compactification of $(U,f)$
and $P$ be the pole divisor of $f$ on $X$.
\begin{enumerate}
\item
The \textit{de Rham cohomology of $(U,\nb)$ with compact support}
is the hypercohomology
\[ \HdRc^i(U,\nb) = \HH^i\left(X, \OO(-S) \xrightarrow{\nb} \OL^1(-S+P)
	\to\cdots\to \OL^p(-S+pP) \to\cdots\right). \]
\item
Write $S = \red{P} + T$.
For any $\lam\in\RR$, define
\[ \cpt{F}^\lam(\nb) := \left(F^\lam(\nb)\right)(-T), \]
regarded as a subcomplex of $F^\lam(\nb)$ on $X$.
(The stability under $\nb$ in $\cpt{F}^\lam(\nb)$
is easy to check.)
Set $\cpt{\Gr}^\lam(\nb) = \cpt{F}^\lam(\nb)/\cpt{F}^{\lam+}(\nb)$.
We let
\begin{eqnarray*}
\cpt{H}^i\left(U,F^\lam(\nb)\right) &:=&
	\HH^i\left(X,\cpt{F}^\lam(\nb)\right) \\
\cpt{H}^i\left(U,\Gr^\lam(\nb)\right) &:=&
	\HH^i\left(X,\cpt{\Gr}^\lam(\nb)\right).
\end{eqnarray*}
\item
By Prop.\ref{Prop:Change-Support}
(for $D=-S, E=\red{P}$),
we have
\[ \HdRc^i(U,\nb) = \cpt{H}^i\left(U,F^0(\nb)\right). \]
The \textit{irregular Hodge filtration} on $\HdRc^i(U,\nb)$
is the filtration
\begin{eqnarray*}
F^\lam \HdRc^i(U,\nb) &=&
	\Image\left\{\HH^i\left(X,\cpt{F}^\lam(\nb)\right) \to
		\HH^i\left(X,\cpt{F}^0(\nb)\right) \right\} \\
	&=& \Image\left\{\cpt{H}^i\left(U,F^\lam(\nb)\right) \to
		\HdRc^i(U,\nb) \right\}.
\end{eqnarray*}
\end{enumerate}
\bigskip

In particular,
if $f:U\to\AA$ is proper (i.e., $S = \red{P}$),
one has the natural isomorphism
\[ \left(\HdRc^i(U,\nb),F^\lam\right) \isoto
	\left(\HdR^i(U,\nb),F^\lam\right). \]

\begin{Prop}\label{Prop:H_c}
Let $U,f,\nb$ be as before.
We have the following functorial properties.
\begin{enumerate}
\item
Let $a: U' \to U$
be a proper morphism of smooth quasi-projective varieties
and let $\nb' = a^*\nb$ be the pullback connection on $U'$.
Then the natural map
$a^*(\Om_U^\bullet,\nb) \to (\Om_{U'}^\bullet,\nb')$
induces
\[ a^*: \cpt{H}^q\left(U,F^\lam(\nb)\right) \to
	\cpt{H}^q\left(U',F^\lam(\nb')\right). \]
\item
Let $i:V\to U$ be a smooth divisor
and $j: U^\circ \to U$ be the complement.
Then we have the natural long exact sequence
\[ \cdots \to \cpt{H}^q\left(U^\circ,F^\lam(\nb)\right) \xrightarrow{j_*}
	\cpt{H}^q\left(U,F^\lam(\nb)\right) \xrightarrow{i^*}
	\cpt{H}^q\left(V,F^\lam(\nb)\right) \to\cdots. \]
\end{enumerate}
\end{Prop}

\pf
(i)
Take a compactification $b: X' \to X$ of $a: U'\to U$
such that $(X',X'\setminus U')$ and $(X,X\setminus U)$
are good compactifications of $(U',f\circ a)$ and $(U,f)$,
respectively.
Write $X'\setminus U' = T'+\red{P'}$
where $P'$ is the pole divisor of $f\circ a$.
Then we have $T' \subset b^{-1}(T)$ since $a$ is proper.
Thus
$b^*\cpt{F}^\lam(\nb) \subset \cpt{F}^\lam(\nb')$
and the assertion follows.

(ii)
Choose a good compactification $(X,S)$ of $(U,f)$
such that $S+\bar{V}$ forms a normal crossing divisor of $X$
where $\bar{V}$ is the closure of $V\subset X$.
Thus $(X,S+\bar{V})$ and $(\bar{V},S\cap\bar{V})$ are good for
$(U^\circ,f)$ and $(V,f)$, respectively.
On $X$, we have
\[ \cpt{F}^\lam(\nb|_{U^\circ}) = \cpt{F}^\lam(\nb)(-\bar{V})
	\subset \cpt{F}^\lam(\nb). \]
Moreover the natural sequence
\[ 0 \to \cpt{F}^\lam(\nb|_{U^\circ}) \to
	\cpt{F}^\lam(\nb) \xrightarrow{i^*}
	\cpt{F}^\lam(\nb|_V) \to 0 \]
is exact as can be derived easily by local computations.
(Cf.~\cite[Example 7.23(1)]{PS} and \cite[Prop.3.7.15]{EZ}
for the case $f$ is trivial.
In both references,
$V$ is allowed to be a normal crossing divisor.)
The assertion then follows by taking hypercohomology.
\qed

\section{The duality}\label{Sect:Duality}

In this section
we assume that $U$ is irreducible of $\dim U = n$.
Recall that,
when we want to emphasis the dependence of $f$,
we write $\nb_f = d+df$
for the twisted connection.
We shall define canonically a perfect bilinear pairing
\[ \HdR^i(U,\nb_f)\times \HdRc^{2n-i}(U,\nb_{-f})
	\xrightarrow{\ggen{\,,\,}} \HdRc^{2n}(U) = \CC \]
for every $i$,
which is compatible with the irregular Hodge filtrations on them.

Let $(X,S)$ be a fixed good compactification of $(U,f)$
throughout the discussion.
As before, let $P$ be the pole divisor of $f$ on $X$ and write
\[ S = \red{P}+T. \]

\subsection{The pairing on the de Rham cohomology}
\newcommand{\cC}{\mathcal{C}}

To define the pairing on the de Rham cohomology of $(U,\nb)$,
we mimic Deligne's construction in \cite[p.124]{DMR}.
\medskip

We construct a chain map
\begin{equation}\label{Eq:dR-pairing-I}
F^0(\nb_f)\tensor{\OO_X} \left(\cpt{F}^0(\nb_{-f})(-\red{P})\right)
	\xrightarrow{\gen{\,,\,}} (\Om_X^\bullet,d)
\end{equation}
to the usual de Rham complex of $X$
in the following way.
First by Prop.\ref{Prop:Change-Support},
the inclusions of complexes
\[ F^0(n) = F^0(\nb_f)(-nP) \hookrightarrow F^0(\nb_f) \]
is a quasi-isomorphism
and we have a chain map
from $F^0(n)\otimes \cpt{F}^0(\nb_{-f})(-\red{P})$:
\begin{equation}\label{Eq:dR-pairing-II}
\begin{array}{c}
\left\{\begin{array}{c}
\OO(-nP) \xrightarrow{\nb_f} \cdots\to
	\OL^{n-i}(-iP) \to\cdots \to \OL^n \\ \bigotimes \\
\OO(-S) \xrightarrow{\nb_{-f}} \cdots\to \OL^j(-S+jP)
	\to\cdots\to \OL^n(-S+nP)
\end{array}\right\} \\ \Big\downarrow \\
\OO(-S-nP) \xrightarrow{d} \cdots\to \OL^{n-k}(-S-kP)
	\to\cdots\to \OL^n(-S)=\Om^n.
\end{array}
\end{equation}
Here the pairings
\[ \OL^{n-i}(-iP) \tensor{\OO_X} \OL^j(-S+jP) \to \OL^{n-i+j}(-S-(i-j)P) \]
appeared in the above chain map
are the natural exterior product.
Now the last complex in \eqref{Eq:dR-pairing-II}
is a subcomplex of $(\Om_X^\bullet,d)$.
Thus, via this inclusion,
we obtain the desired chain map \eqref{Eq:dR-pairing-I}.
\medskip

Taking hypercohomology,
we then obtain the Poincar\'e pairing
\[ \HdR^i(U,\nb_f) \times \HdRc^{2n-i}(U,\nb_{-f}) \xrightarrow{\ggen{\,,\,}}
	\HdR^{2n}(X) = H^n(X,\Om^n). \]

\begin{Thm}\label{Thm:Poincare-Pairing}
For any $i$,
the Poincar\'e pairing $\ggen{\,,}$ constructed above is perfect.
\end{Thm}

\pf
Indeed we have the perfect pairing
\[ \OL^i(mS) \otimes \OL^{n-i}(m'S) \to \OL^n((m+m')S)
	= \Om^n_X((m+m'+1)S). \]
Consider the Hom-sheaf with value in $\Om^n_X$
\[ (\bullet)^\wedge := \Homsh_{\OO_X}(\bullet, \Om^n_X). \]
Then we have
\begin{eqnarray*}
\left(\cpt{F}^0(\nb_{-f})\right)^\wedge &=&
	\left[\OO \xrightarrow{\nb_f} \cdots \to \posz{\OL^n(nP)}\right]
	\otimes\OO(-nP) \\
&\simeq& F^0(\nb_f)[n]
	\quad(\text{Prop.\ref{Prop:Change-Support}}).
\end{eqnarray*}
Therefore by filtering the complexes and the Serre duality,
we have
\[ \HdRc^i(U,\nb_{-f})^\vee = \HH^{n-i}\left(X,F^0(\nb_f)[n]\right)
	= \HdR^{2n-i}(U,\nb_f). \]
Here $(\bullet)^\vee$ denotes the dual vector space.

One can use,
e.g., the fine resolution of the twisted de Rham complex
into sheaves of $\cC^\infty$ $(p,q)$-forms
(with appropriate poles along $S$)
to check that
the argument here
is compatible with the definition of the pairing.
(Cf.~the proof of the next theorem.)
{\qed}

\begin{Thm}
The two pairs
$(\HdR^i(U,\nb_f),F^\bullet)$ and $(\HdRc^{2n-i}(U,\nb_{-f}),F^\bullet)$
of filtered vector spaces are dual to each other
via the perfect Poincar\'e pairing
(up to a degree shift).
More precisely, for any $\lam$ we have
\begin{equation}\label{Eq:Pairing-Fil}
\ggen{F^\lam\HdR^i(\nb_f), F^{(n-\lam)+}\HdRc^{2n-i}(\nb_{-f})} = 0
	= \ggen{F^{\lam+}\HdR^i(\nb_f), F^{(n-\lam)}\HdRc^{2n-i}(\nb_{-f})}
\end{equation}
and the Poincar\'e pairing induces a duality between
$\Gr^\lam\HdR^i(\nb_f)$ and $\Gr^{n-\lam}\HdRc^{2n-i}(\nb_{-f})$.
(We have omitted the base $U$ inside the cohomology
in the formulas.)
\end{Thm}

\pf
We use the fine resolution into $\cC^\infty$ $(p,q)$-forms.
Suppose $\om \in F^\lam \HdR^i(\nb_f)$.
Since $\om$ is the image of an element in $H^i(F^\lam(\nb_f))$,
it is represented by
\[ \sum \om_p \in \bigoplus_{p\geq\lam}
	\Gamma\left(X,\OL_\infty^{p,i-p}(\flr{(p-\lam)P})\right)
	\subset \bigoplus_{p\geq 0}
	\Gamma\left(X,\OL_\infty^{p,i-p}(\flr{(p-\lam)P})\right). \]
Here $\OL_\infty^{p,i-p}$ denotes
the sheaf of logarithmic $(p,i-p)$-forms
with $\cC^\infty$ coefficients.
As the inclution
\[ \left(\OL^p(-\red{P}+\flr{(p-\lam)P}),\nb_f\right)_{p\geq 0} \subset
	\left(\OL^p(\flr{(p-\lam)P}),\nb_f\right)_{p\geq 0} \]
is a quasi-isomorphism,
there exists $\alpha \in \bigoplus_{p\geq 0}
	\Gamma\left(X,\OL_\infty^{p,i-1-p}(\flr{(p-\lam)P})\right)$
such that
\[ \sum(\om+D(\alpha))_p \in \bigoplus_{p\geq 0}
	\Gamma\left(X,\OL_\infty^{p,i-p}(-\red{P}+\flr{(p-\lam)P})\right). \]
Here $D$ is the total differential given by
\[ D = \nb + (-1)^p\bar{d}
\quad\text{on $\OL_\infty^{pq}(rP)$}. \]

Now given $\eta \in F^{(n-\lam)+} \cpt{H}^{2n-i}(\nb_{-f})$,
which is represented as the sum
\[ \sum \eta_q \in \bigoplus_{q>n-\lam}
	\Gamma\left(X,\OL_\infty^{q,2n-i-q}(-T+\flr{(q-n+\lam-\eps)P})\right)
	\quad\text{for some $\eps >0$}, \]
one obtains
\begin{eqnarray*}
\ggen{\om,\eta} &=& \sum_{q>n-\lam} \ggen{(\om+D(\alpha))_{n-q},\eta_q}
		\quad(\text{definition}) \\
	&=& \sum_{q>n-\lam} \ggen{D(\alpha)_{n-q},\eta_q}
		\quad(\text{as $\om_{n-q}=0$ if $n-q < \lam$}) \\
	&=& \int_X D(\alpha\wedge\eta)
		\quad(\text{as $D\eta = 0$}) \\
	&=& 0
		\quad(\text{Stokes}).
\end{eqnarray*}
The application of the Stokes theorem in the last equality above is valid
since the $(2n-1)$-form $\alpha\wedge\eta$
has no poles on the compact $X$.
Indeed for $q>n-\lam$ and any $\eps >0$,
\begin{eqnarray*}
\alpha_{n-q}\wedge\eta_q &\in&
	\Gamma\left(X,
	\OL_\infty^{n,n-1}(\flr{(n-q-\lam)P}-T+\flr{(q-n+\lam-\eps)P})\right) \\
	&\subset& \Gamma\left(X,\OL_\infty^{n,n-1}(-\red{P}-T)\right) \\
	&=& \Gamma\left(X,\Om_\infty^{n,n-1}\right),
\end{eqnarray*}
while similarly
\begin{eqnarray*}
\alpha_{n-1-q}\wedge\eta_q &\in&
	\Gamma\left(X,
	\OL_\infty^{n-1,n}(\flr{(n-1-q-\lam)P}-T+\flr{(q-n+\lam-\eps)P})\right) \\
	&\subset& \Gamma\left(X,\OL_\infty^{n-1,n}(-\red{P}-T)\right) \\
	&\subset& \Gamma\left(X,\Om_\infty^{n-1,n}\right).
\end{eqnarray*}

The second equality of \eqref{Eq:Pairing-Fil}
can be proved similarly.
\qed

\begin{Cor}
The filtered cohomology $(\HdRc^i(U,\nb),F^\lam)$ with compact support
does not dependent on the choice of the compactification $X$.
\end{Cor}

\subsection{Pairings on the spectral sequence}
Take a sequence
\[ \lam_{-1}<0=\lam_0<\lam_1<\cdots<\lam_N=n<\lam_{N+1}
\quad\text{with $\lam_i+\lam_{N-i} = n$} \]
where $\lam_0,\cdots,\lam_N$ are
all the non-negative jumps of the filtration $F^\lam(\nb)$ on $X$.
Notice that we have
\[ \flr{\lam_{i-1}P} + \flr{-\lam_iP} = -\red{P} \]
as can be checked easily.
The associated Hodge to de Rham spectral sequence reads
\[ E_{1,*}^{p,q} = H_*^{p+q}\left(U,\Gr^{\lam_p}(\nb)\right)
	\Longrightarrow H_\mathrm{dR,*}^{p+q}(U,\nb) \]
where $*=\mathrm{c}$ or nothing.

\subsubsection*{$E_r$-terms and jumping gradings}
For $0\leq i<j\leq N+1$, let
\[ G\binom{i}{j} := \tot
	\left[\begin{array}{c} F^{\lam_j}(\nb_f) \\ \downarrow \\
	\posz{F^{\lam_i}(\nb_f)} \end{array}\right]
		\quad\text{and}\quad
\cpt{G}\binom{i}{j} := \tot
	\left[\begin{array}{c} \cpt{F}^{\lam_j}(\nb_{-f}) \\ \downarrow \\
	\posz{\cpt{F}^{\lam_i}(\nb_{-f})} \end{array}\right] \]
be the \textit{jumping gradings},
which are the representatives of the quotients
$F^{\lam_i}(\nb_f)/F^{\lam_j}(\nb_f)$ and
$\cpt{F}^{\lam_i}(\nb_{-f})/\cpt{F}^{\lam_j}(\nb_{-f})$,
respectively.
Notice that similar to Cor.\ref{Cor:HGr-well-defined}
the hypercohomology of $G\binom{i}{j}$ over $X$
is independent of the choice of $X$.

For two pairs of numbers $i<j$ and $i'<j'$,
we say $(i,j)\geq(i',j')$
if $i\geq i'$ and $j\geq j'$.
Then for any $(i,j)\geq(i',j')$,
there is the natural componentwise inclusion
\begin{equation}\label{Eq:G-G}
G\binom{i}{j} \to G\binom{i'}{j'}
\end{equation}
and if $i<j<k$, one has the distinguished triangle
\[ G\binom{j}{k} \to G\binom{i}{k} \to G\binom{i}{j} \xrightarrow{+1}. \]
We have
\[ E_r^{p,q-p} = \Image\left\{
	\HH^q\left(X,G\binom{p}{p+r}\right) \to
	\HH^q\left(X,G\binom{p-r+1}{p+1}\right) \right\}, \]
traditionally regarded as a subquotient of
$\HH^q\left(X,G\binom{p}{p+1}\right)$.
The following commutative diagram
illustrates the various terms in the spectral sequence.
\begin{equation}\label{Eq:Er-in-SS}
\xymatrix{
\HH^q(G\binom{p}{N+1}) \ar[r]
	\ar@/_2pc/[rrrrddd]_{\text{Image = $E_\infty$}} &
\HH^q(G\binom{p}{N}) \ar[r] & \cdots \ar[r] &
\HH^q(G\binom{p}{p+2}) \ar[r]\ar@/_1pc/[rd]_(.3){\text{Image = $E_2$}} &
\HH^q(G\binom{p}{p+1}) = E_1^{p,q-p} \ar[d] \\
&&&& \HH^q(G\binom{p-1}{p+1}) \ar[d] \\
&&&& \vdots \ar[d] \\
&&&& \HH^q(G\binom{0}{p+1}) }
\end{equation}
(Here we omit the base $X$ of the hypercohomology.
Notice that $G\binom{p}{N+1} = F^{\lam_p}(\nb_f)$
and $G\binom{p}{p+1} = \Gr^{\lam_p}(\nb_f)$.)
\medskip

Similar pictures hold for $\cpt{G}\binom{i}{j}$
(but replace $\nb_f$ by $\nb_{-f}$ in this case).

\subsubsection*{From subs to quots}
Recall the complex $F^0(\lam)$ defined
in \eqref{Eq:F0lam} of \S\ref{Sect:IrrHodge}.
Notice that by Prop.\ref{Prop:Change-Support} the inclusion
\[ F^0(\lam) = F^0(\nb)(\flr{-\lam P}) \to F^0(\nb) \]
is a quasi-isomorphism for any $\lam\geq 0$.
Define
\[ Q_\lam = F^0(\lam)^{<\ceil{\lam}}. \]
We have the short exact sequence
\[ 0 \to F^\lam(\nb_f) \to F^0(\lam) \to Q_\lam \to 0. \]

Now for $0\leq i<j\leq N+1$,
define
\[ Q\binom{i}{j} := \tot
	\left[\begin{array}{c} Q_{\lam_j} \\ \downarrow \\
	\posz{Q_{\lam_i}} \end{array}\right]. \]
Then for $(i,j) \geq (i',j')$, we have
the componentwise quotient
\begin{equation}\label{Eq:Q-Q}
Q\binom{i}{j} \to Q\binom{i'}{j'}
\end{equation}
and the distinguished triangles
\[ G\binom{i}{j} \to
\tot\left[\begin{array}{c} F^0(\lam_j) \\ \downarrow \\
	\posz{F^0(\lam_i)} \end{array}\right]
\to Q\binom{i}{j} \xrightarrow{+1}. \]
Since the middle term is quasi-isomorphic to zero,
we obtain a quasi-isomorphism
\begin{equation}\label{Eq:QI-Q-G}
Q\binom{i}{j}[-1] \isoto G\binom{i}{j}.
\end{equation}
The natural maps \eqref{Eq:G-G} and \eqref{Eq:Q-Q}
are compatible with the above quasi-isomorphism.

\subsubsection*{The pairings}
To show that the spectral sequence is compatible with the duality,
we should construct pairings
\begin{equation}\label{Eq:Pairing-i-j}
\gen{\,,}^i_j: G\binom{i}{j} \tensor{\OO_X}
	\cpt{G}\binom{N+1-j}{N+1-i}
\to (\Om_X^\bullet,d)
\end{equation}
for all $0\leq i<j\leq N+1$
which induce perfect pairings on cohomology
and are compatible with respect to the partial ordering of various $(i,j)$.

Define a pairing
\[ Q\binom{i}{j}[-1] \tensor{\OO_X} \cpt{G}\binom{N+1-j}{N+1-i}
	\xrightarrow{\gen{\,,\,}^i_j} (\Om_X^\bullet,d) \]
as follows.
For $\om$ in the degree $p$ term of $Q\binom{i}{j}[-1]$
\begin{equation}\label{Eq:p-in-G}
\om = (\om_1,\om_2) \in
	\OL^p\big(\flr{(p-\lam_j)P}\big)\oplus
		\OL^{p-1}\big(\flr{(p-1-\lam_i)P}\big)
\end{equation}
and $\eta$ in the degree $q$ term of $\cpt{G}\binom{N+1-j}{N+1-i}$
\begin{equation}\label{Eq:q-in-Gc}
\begin{split}
\eta = (\eta_1&,\eta_2) \\
	&\in \OL^{q+1}\big(-T+\flr{(q+1-\lam_{N+1-i})P}\big)\oplus
		\OL^q\big(-T+\flr{(q-\lam_{N+1-j})P}\big) \\ 
	&= \OL^{q+1}\big(-T+\flr{(q+1-n+\lam_{i-1})P}\big)\oplus
		\OL^q\big(-T+\flr{(q-n+\lam_{j-1})P}\big),
\end{split}\end{equation}
we set
\begin{eqnarray*}
\gen{\om,\eta}^i_j := \om_1\eta_2+(-1)^q\om_2\eta_1 &\in&
	\OL^{p+q}\big(-T+\flr{\lam_{i-1}P}+\flr{-\lam_iP}+(p+q-n)P\big) \\
	&=& \OL^{p+q}\big(-S-(n-p-q)P\big) \\
	&\subset& \Om_X^{p+q}.
\end{eqnarray*}
One checks readily that
\[ d\gen{\om,\eta}^i_j = \gen{\nb_f(\om),\eta}^i_j
	+ (-1)^p\gen{\om,\nb_{-f}(\eta)}^i_j, \]
where we have adapted the sign convention
\begin{eqnarray*}
\nb_f(\om_1,\om_2) &=&
	\left(\nb_f(\om_1),(-1)^{p-1}\om_1+\nb_f(\om_2)\right) \\
\nb_{-f}(\eta_1,\eta_2) &=&
	\left(\nb_{-f}(\eta_1),(-1)^q\eta_1+\nb_{-f}(\eta_2)\right).
\end{eqnarray*}
(The sign $(-1)^{p-1}$ in the first equation
is due to the shift $[-1]$ in $Q\binom{i}{j}[-1]$.)

Using the quasi-isomorphism \eqref{Eq:QI-Q-G}
we obtain the pairing \eqref{Eq:Pairing-i-j}
which induces a pairing
\begin{equation}\label{Eq:Pairing-H-i-j}
\HH^q\left(X,G\binom{i}{j}\right) \times
	\HH^{2n-q}\left(X,\cpt{G}\binom{N+1-j}{N+1-i}\right)
	\xrightarrow{\ggen{\,,\,}^i_j} \HdR^{2n}(X) = \CC
\end{equation}
for each $q$.

\begin{Thm}
For all $0\leq i<j\leq N+1$,
the parings $\ggen{\,,}^i_j$ are perfect
and they are compatible with each other
under the partial ordering of $(i,j)$ and the map \eqref{Eq:G-G}.
\end{Thm}

\pf
The proof of the perfectness
is similar to that of Thm.\ref{Thm:Poincare-Pairing}.
One shows by induction on the length $l$ that
the cohomology of the various truncations
\[ Q\binom{i}{j}[-1]^{\geq n-l}
\quad\text{and}\quad
\cpt{G}\binom{N+1-j}{N+1-i}^{\leq l} \]
are dual to each other via the pairing.
In each step,
the perfectness follows from the classical Serre duality
asserting the perfectness of the pairing
\[ H^q\left(X,\OL^p(D)\right) \times H^{n-q}\left(X,\OL^{n-p}(-S-D)\right)
	\to H^n(X,\Om^n) = \HdR^{2n}(X). \]

The compatibilities with the ordering of $(i,j)$
and with the definition of the pairings are clear
by e.g.~writing everything in terms of $\cC^\infty$ differential forms.
\qed

\bigskip
\noin\textit{Remark.}
For $(i,j) = (0,N+1)$,
the quasi-isomorphism \eqref{Eq:QI-Q-G}
reduces to the inclusion $F^0(\lam_{N+1}) \to F^0(\nb)$.
One can use this inclusion
and the complex $\cpt{F}^0(\nb_{-f})$
instead of $F^0(n) \to F^0(\nb)$ and $\cpt{F}^0(\nb_{-f})(-\red{P})$,
respectively
in the construction \eqref{Eq:dR-pairing-II}
of the previous subsection
to define the (same) pairing on the de Rham cohomology.
In this case the above theorem
then recovers Thm.\ref{Thm:Poincare-Pairing}.

\begin{Cor}
The cohomology $\cpt{H}^q\left(U,F^\lam(\nb_f)\right)$
and $\cpt{H}^q\left(U,\Gr^\lam(\nb_f)\right)$
do not depend on the choice of the good compactification
$(X,S)$ of $(U,f)$.
\end{Cor}

\pf
Applying the above theorem for $i=0$,
we see that $\cpt{H}^{2n-q}\left(U,F^{\lam_{N+1-j}}(\nb_{-f})\right)$
is canonically dual to $\HH^q\left(X,G\binom{0}{j}\right)$.
As already mentioned,
this later space is independent of the choice of $X$.
Thus after renaming the indices and the function $f$,
we see that $\cpt{H}^q\left(U,F^\lam(\nb_f)\right)$
is independent of $X$.
The other statement follows by taking long exact sequence
and from the compatibility of the pairings
with respect to the partial ordering of $(i,j)$.
\qed
\bigskip

Using the description in \eqref{Eq:Er-in-SS},
the above theorem and the remark after it imply the following.

\begin{Cor}
The pairings $\ggen{\,,\,}^i_j$ induce perfect pairings
\[ \ggen{\,,}_r : E_r^{p,q} \times E_{r,\mathrm{c}}^{N-p,2n-N-q} \to \CC. \]
They are compatible with the pairing on the de Rham cohomology
in the sense that
\[ \ggen{\om,\eta}_r = \ggen{\om,\eta} \]
for $(\om,\eta) \in H^{p+q}\left(U,F^{\lam_p}(\nb)\right) \times
	\cpt{H}^{2n-p-q}\left(U,F^{n-\lam_p}(\nb)\right)$
projected into the $E_r$-term in the left
and into the de Rham cohomology in the right.
\end{Cor}

\begin{Prop}
For $0\leq i<j<k\leq N+1$,
the sequence of pairings constructed in \eqref{Eq:Pairing-H-i-j}
on the cohomology of the two distinguished triangles
\[\xymatrix{
& G\binom{j}{k} \ar[r] & G\binom{i}{k} \ar[r] &
	G\binom{i}{j} \ar[r]^{+1} & \\
& \ar[l]_(.7){+1} \cpt{G}\binom{N+1-k}{N+1-j} & \ar[l]
	\cpt{G}\binom{N+1-k}{N+1-i} & \ar[l] \cpt{G}\binom{N+1-j}{N+1-i} & }\]
commutes up to sign.
\end{Prop}

\pf
We show that already in the chain level
the diagram
\[\xymatrix{
G\binom{i}{j}\otimes\cpt{G}\binom{N+1-j}{N+1-i}
	\ar[rr]^(.6){\gen{\,,\,}^i_j}\ar@<-7ex>[d]_{\delta} &&
		(\Om_X^\bullet,d) \ar@{=}[d] \\
G\binom{j}{k}\otimes\cpt{G}\binom{N+1-k}{N+1-j}
	\ar[rr]^(.6){\gen{\,,\,}^j_k}\ar@<-4ex>[u]_{\cpt{\delta}} &&
		(\Om_X^\bullet,d) }\]
commutes up to sign
where $\delta$ and $\cpt{\delta}$ are the map $+1$
in the distinguished triangles.
Let $\om$ and $\eta$
be degree $p$ and $q$ elements
in $Q\binom{i}{j}[-1]$ and $\cpt{G}\binom{N+1-k}{N+1-j}$
as in \eqref{Eq:p-in-G} and \eqref{Eq:q-in-Gc},
respectively
(but replace $(i,j)$ by $(j,k)$ in \eqref{Eq:q-in-Gc}).
Then
\[\begin{split}
\delta(\om) = (&0,\om_1) \\ &\in
	\OL^{p+1}\big(\flr{(p+1-\lam_k)P}\big)\oplus
		\OL^p\big(\flr{(p-\lam_j)P}\big), \\
\cpt{\delta}(\eta) = (&0,\eta_1) \\ &\in
	\OL^{q+2}\big(-T+\flr{(q+2-n+\lam_{i-1})P}\big)\oplus
		\OL^{q+1}\big(-T+\flr{(q+1-n+\lam_{j-1})P}\big).
\end{split}\]
Thus
\[ \gen{\delta(\om),\eta}^j_k = (-1)^q\om_1\eta_1
	= (-1)^q\gen{\om,\cpt{\delta}(\eta)}^i_j. \]

On the cohomology level,
we then have
\[ \ggen{\delta(\om),\eta}^j_k + (-1)^p\ggen{\om,\cpt{\delta}(\eta)}^i_j = 0 \]
for $\om \in \HH^p\left(X,G\binom{i}{j}\right)$
and $\eta \in \HH^{2n-p-1}\left(X,\binom{N+1-k}{N+1-j}\right)$.
\qed
\bigskip

\noin\textit{Remark.}
The above proposition
reduces to the special case of the compatibility
of the pairings with respect to the ordering of $(i,j)$
in the previous theorem
if the spectral sequence degenerates at the $E_1$-terms,
since the connection maps in cohomology
induced by $\delta$ and $\cpt{\delta}$ in the proof
are then all zero.

\section{The hypersurface case}\label{Sect:Hypersurf}

We first recall the following well-known relation
between the exponential sums and counting solutions of equations
over a finite field.
Let $f$ be a regular function on a quasi-projective variety $U$
over a finite field $\kappa$.
To count the number $N(f)$ of elements of the zero set
\[ \{ x \in U(\kappa) \,|\, f(x)=0 \}, \]
one brings in a non-trivial additive character
$\chi: \kappa \to \CC^\times$
and introduces a new variable $z=$
a fixed coordinate of an affine line $\AA$ over $\kappa$.
Let $\wt{f}=zf$, a regular function on $\AA\times U$.
Then we have
\[ \sum_{x\in (\AA\times U)(\kappa)} \chi\left(\wt{f}(x)\right)
	= q\cdot N(f) \]
where $q$ is the cardinality of $\kappa = \AA(\kappa)$.

Now the exponential sum in the left hand term of the equality above
is related to the finite-field counterpart of the twisted de Rham cohomology
while the right hand term consists of
information of the closed subscheme of $U$ defined by $f$.
This suggests that in the world over the field of complex numbers,
the de Rham cohomology of the connection $\nb_{\wt{f}}$
over the product $\AA\times U$
together with its irregular Hodge filtration
should reflect the usual de Rham cohomology
of the closed subscheme $(f)$ defined by $f$
with the usual Hodge filtration.
We work out this analogue in this section.
We consider the case
where $\red{f}$ defines a smooth divisor of $U$
since we only have defined the twisted de Rham complexes
and the filtrations
for smooth varieties.

\newcommand{\pr}{\operatorname{pr}}

\begin{Lemma}\label{Lemma:SimpleKunneth}
Let $U$ be quasi-projective and smooth
and $\wt{U} = \AA\times U$.
Consider the two projections
\[ \AA \xleftarrow{a} \wt{U} \xrightarrow{b} U. \]
\begin{enumerate}
\item
Let $\nb$ be the twisted connection on $\AA$
associated with the identity map.
Let
$\wt{\nb} = a^*\nb = \nb\boxtimes d$
be the pullback connection on $\wt{U}$.
Then for any $i,\lam$ we have
\[ H^i\left(\wt{U},F^\lam(\wt{\nb})\right) = 0. \]
\item
Let $\nb$ be the twisted connection associated with
a regular function on $U$
and $\wt{\nb} = b^*\nb = d\boxtimes\nb$
be the pullback.
Then for any $i,\lam$ we have
\[ H^i\left(\wt{U},F^\lam(\wt{\nb})\right) = H^i\left(U,F^\lam(\nb)\right). \]
\end{enumerate}
\end{Lemma}

\pf
Let $(X,S)$ be a good compactification of $(U,f)$
where $f=0$ in (i) or the regular function in (ii).
Then $(\PP\times X,\{\infty\}\times X \cup \PP\times S)$
is in fact a good compactification of
$(\wt{U},\mathrm{id}\circ a)$ in (i)
or $(\wt{U},f\circ b)$ in (ii).

(i) Let $F_\boxtimes^\lam$
be the exterior product filtration on $\PP\times X$
of $F^\mu(\nb)$ on $\PP$
and $F^\nu(d)$ on $X$.
By Prop.\ref{Prop:Prod-Fil}
the natural inclusion
$F_\boxtimes^\lam \to F^\lam(\wt{\nb})$ is a quasi-isomorphism.
Thus we only need to compute the hypercohomology of
$F_\boxtimes^\lam$.
On the other hand,
on the good compactification $\PP$ of $\AA$
we have
\[ \Gr^0(\nb) \cong \OO_\PP(-1)
\quad\text{and}\quad
F^{0+}(\nb) = \OL^1_{\AA\subset\PP}[-1] \cong \OO_\PP(-1)[-1]. \]
Thus $F_\boxtimes^\lam$ is quasi-isomorphic to an extension of
\[ A:=
\Gr^0(\nb)\boxtimes F^\lam(d) \cong \OO_\PP(-1)\boxtimes F^\lam(d) \]
by
\[ B:= F^{0+}(\nb)\boxtimes F^{\lam-}(d) \cong
	\OO_\PP(-1)[-1]\boxtimes F^\lam(d). \]
Since $A$ and $B$ have trivial hypercohomology,
the assertion follows.

(ii) Similarly let $F_\boxtimes^\lam$
be the product filtration on $\PP\times X$
of $F^\mu(d)$ on $\PP$
and $F^\nu(\nb)$ on $X$.
We have the quasi-isomorphism
$F_\boxtimes^\lam \isoto F^\lam(\wt{\nb})$
by Prop.\ref{Prop:Prod-Fil}.
This time on $\PP$ we have
\[ \Gr^0(d) \cong \OO_\PP
\quad\text{and}\quad
F^{0+}(d) = \OL^1_{\AA\subset\PP}[-1] \cong \OO_\PP(-1)[-1]. \]
Thus $F_\boxtimes^\lam$ on $\PP\times X$
has the same hypercohomology
as $F^\lam(\nb)$ on $X$.
\qed

\begin{Lemma}\label{Lemma:on_Ucirc}
Let $f$ be a nowhere vanishing regular function
on a smooth quasi-projective $U^\circ$
and $\wt{f}=zf$ on $\AA\times U^\circ$ where $z=$ identity on $\AA$.
Let $\nb = d + d\wt{f}$ be the twisted connection on $\AA\times U^\circ$.
Then for all $i,\lam$ we have
\[ H^i\left(\AA\times U^\circ,F^\lam(\nb)\right) = 0. \]
\end{Lemma}

\pf
We have the commutative diagram
\[\xymatrix{
\AA\times U^\circ \ar[r]^\alpha\ar[dr]_{\wt{f}=zf} &
	\AA\times U^\circ \ar[d]^{\pr_1} \\
& \AA  }\]
where
\[ \alpha(z,x) := (zf(x), x), \]
is an isomorphism.
Thus to prove the assertion,
one reduces to the case where $\wt{f} = z$
via the isomorphism $\alpha$.
The assertion then follows from
Lemma \ref{Lemma:SimpleKunneth}(i).
\qed

\begin{Thm}
Consider a pair $(U,f)$ as before.
Let $V = \red{f}$ be the closed subvariety of $U$ defined by $f$.
Let $\wt{f} = zf$ on $\AA\times U$ where $z=$ identity on $\AA$.
Assume that $V$ is smooth.
Then, for any $i,\lam$,
\[ H_\mathrm{dR,*}^{i+2}\left(\AA\times U,\nb_{\wt{f}}\right) =
	H_\mathrm{dR,*}^i(V) \]
and
\[ H_*^{i+2}\left(\AA\times U,\Gr^{\lam+1}(\nb_{\wt{f}})\right) =
	H_*^{i-\ceil{\lam}}\left(V, \Om^{\ceil{\lam}}\right) \]
where $*=\mathrm{c}$ or nothing.
\end{Thm}

\pf
By duality,
it is enough to consider the case for the cohomology with compact support.

Let $U^\circ = U\setminus V$.
Thus the three
\[ \AA\times U^\circ \hookrightarrow
	\AA\times U \hookleftarrow \AA\times V \]
form an open-closed decomposition
and, by Prop.\ref{Prop:H_c}(ii), we have the long exact sequence
\[ \cdots\to \cpt{H}^i\left(\AA\times U^\circ,F^\lam(\nb)\right) \to
	\cpt{H}^i\left(\AA\times U,F^\lam(\nb)\right) \to
	\cpt{H}^i\left(\AA\times V,F^\lam(d)\right) \to\cdots. \]
By the dual of the lemma above and the K\"unneth formula,
we then have
\[ \cpt{H}^i\left(\AA\times U,F^\lam(\nb)\right) =
	\cpt{H}^i\left(\AA\times V,F^\lam(d)\right) =
	\cpt{H}^{i-2}\left(V,F^{\ceil{\lam}-1}(d)\right). \]
The assertions now follow.
\qed\bigskip

\noin\textit{Remark.}
The above theorem implies in particular that
the Hodge to de Rham spectral sequence degenerates
at $E_1$-terms in the case $\nb = \nb_{\wt{f}}$.
The fact that the filtration $\cpt{F}^\lam(d)$
indeed induces the Hodge filtration
of the canonical mixed Hodge structure on $\HdRc(V)$
can be found in
\cite[\S 4.3.3 and Prop.4.3.6]{EZ}.

\bigskip
Using the same idea,
one has the following statements,
also motivated by counting the number of the solutions
of equations over finite fields.

\begin{Cor}
Let $f_1,\cdots,f_n$ be regular functions on $U$
and consider $\wt{f} = \sum_{i=1}^n z_if_i$ on $\AA^n\times U$
where $\{z_i\}=$ the cartesian coordinates on $\AA^n$.
Suppose that $\sum_{i=1}^n\red{f_i}$ is a
strict normal crossing divisor.
Let $W = \bigcap_{i=1}^n \red{f_i}$.
Then, for any $j,\lam$,
\[ H^{j+2n}_{\mathrm{dR},*}\left(\AA^n\times U,\nb_{\wt{f}}\right) =
	H^j_{\mathrm{dR},*}(W) \]
and
\[ H^{j+2n}_*\left(\AA^n\times U,\Gr^{\lam+n}(\nb_{\wt{f}})\right) =
	H^{j-\ceil{\lam}}_*\left(W, \Om^{\ceil{\lam}}\right) \]
where $*=\mathrm{c}$ or nothing.
\end{Cor}

\pf
Again by duality,
it is enough to consider $*=\mathrm{c}$.
The above theorem gives the results for $n=1$.

In general,
one considers the commutative diagram
\[\xymatrix{
\AA^n\times V \ar[r] & \AA^n\times U &
\AA^n\times U^\circ \ar[l]\ar[d]_\alpha\ar[dr]^{\wt{f}|_{\AA^n\times U^\circ}} \\
&& \AA^n\times U^\circ \ar[r]_(.6){\pr_1} & \AA }\]
where $V = \red{f_1}$, $U^\circ = U\setminus V$ and
\[ \alpha(z_1,\cdots,z_n;x) :=
	(z_1f_1(x)+\cdots+z_nf_n(x),z_2,\cdots,z_n;x) \]
defines an isomorphism.
Now the triangle in the diagram shows that
on $\AA^n\times U^\circ$,
the connection is isomorphic to $d+dz_1$.
By Lemma \ref{Lemma:on_Ucirc}
(applied to $f=1$ on $\AA^{n-1}\times U^\circ$),
the cohomology of this connection
and of its filtered pieces
all vanish.
Therefore the long exact sequence associated
with the open-closed decomposition
in the upper row of the diagram gives
\[ \HdRc^{j+2}\left(\AA^n\times U,F^\lam(\nb)\right) =
	\HdRc^{j+2}\left(\AA^n\times V,F^\lam(\nb|_{\AA^n\times V})\right) =
	\HdRc^j\left(\AA^{n-1}\times V,F^{\lam-1}(\nb')\right) \]
for any $j,\lam$ with
$\nb' = d+d(\sum_{i=2}^nz_if_i)$ on $\AA^{n-1}\times V$.
Here the second equality follows from the fact that
$\nb|_{\AA^n\times V} = d_{z_1}\boxtimes\nb'$
and by the dual of Lemma \ref{Lemma:SimpleKunneth}(ii).
The statements now follow by induction
on the number $n$ of the defining equations of $W$.
\qed

\section{The toric case}\label{Sect:Toric}
\newcommand{\NPF}{F_{\mathrm{NP}}}
\newcommand{\GrNP}{\Gr_{\mathrm{NP}}}
\newcommand{\torX}{X_{\mathrm{tor}}}

Suppose $U$ is a torus.
Inspired by the investigation \cite{AS-Annals}
of exponential sums over a torus
via Dwork's $p$-adic methods
and the work of Kouchnirenko \cite{K} on the Milnor numbers
of isolated singularities,
Adolphson and Sperber in \cite{AS-Nagoya}
study the twisted de Rham cohomology on $U$
(in fact in a more general setting
which also allows multiplicative twists).
They derive that for generic $f$,
the twisted de Rham cohomology is concentrated in a single degree.
The method there is to introduce a filtration,
already appeared in \cite{K},
on the de Rham chain complex
and show that the associated graded complex
has non-trivial cohomology only at one degree.
In this section
we recall their filtration
and show that the induced filtration on the de Rham cohomology
coincides with our irregular Hodge filtration
for $f$ generic.

Our reference for the theory of toric varieties is \cite{F}.
In particular,
see \cite[p.48]{F}
for the existence of the equivariant resolution of singularities
and \cite[p.61]{F}
for the computation of the valuation of a function
on a toric divisor.
\medskip

In this section we let
\[ U = (\AA\setminus 0)^n
\quad\text{and}\quad
f\in\OO(U) = \CC[x^{\pm 1}] \]
where $x=(x_1,\cdots,x_n)$
is the system of cartesian coordinates of $U$.
Recall the following.

\bigskip
\noin\textit{Definition.}
Write $f = \sum_{\alpha \in \ZZ^n} c(\alpha)x^\alpha$.
\begin{enumerate}
\item The \textit{Newton polyhedron} $\Delta(f)$ of $f$
is the convex hull in $\RR^n$ of the finite set
\[ \{0\} \cup \{ \alpha\in\ZZ^n \,|\, c(\alpha) \neq 0 \}. \]
\item The function $f$ is called
\textit{non-degenerate with respect to $\Delta(f)$}
if for any face $\delta$ of $\Delta(f)$
with $0\not\in\delta$,
the system of equations
\begin{equation}\label{Eq:Non-deg-cond}
f_\delta = \frac{\pt f_\delta}{\pt x_1} = \cdots =
	\frac{\pt f_\delta}{\pt x_n} = 0
\end{equation}
has no solution on $U$
where $f_\delta := \sum_{\alpha\in\delta} c(\alpha)x^\alpha$.
\end{enumerate}
\bigskip

One regards $\Delta(f)$ as sitting in the space of characters
$M_\RR := \RR\tensor{\ZZ}\mathrm{Hom}(U,\CC^\times)$.
It then defines a fan on the dual space
$N_\RR := \mathrm{Hom}_\RR(M,\RR)$
where each codimension one face of $\Delta(f)$
corresponds to a ray in the fan,
pointing to the inward normal direction
with respect to the natural pairing
$N_\RR\times M_\RR \to \RR$.
Now one can refine and enlarge the fan
to make a cone decomposition of $N_\RR$
such that the associated toric variety $\torX$
is smooth and proper
and the toric boundary $S:=\torX\setminus U$
is a simple normal crossing divisor of $\torX$.
Each ray in this refined fan
corresponds to an irreducible component of $S$.
We fix this $\torX$ in the sequel.
We have the commutative diagram
\[\xymatrix{
U \ar[r]^f\ar[d] & \AA \ar[d] \\
\torX \ar@{-->}[r] & \PP }\]
where the two vertical arrows are the inclusions
but the lower arrow is just a rational function on $\torX$ in general.

The connection $\nb$ on $U$ again extends
to the complex on $\torX$
\[ (\Om^\bullet(*S),\nb) = \left[
	\OO(*S) \xrightarrow{\nb} \Om^1(*S) \to\cdots\to
		\Om^n(*S)\right] \]
and we have
\[ \HdR^i(U,\nb) = \HH^i\big(\torX,(\Om^\bullet(*S),\nb)\big). \]
\medskip

If $\dim\Delta(f) <n$,
there is a decomposition $U = U'\times U''$
of $U$ into two tori and $f'\in\OO(U')$
such that $f=f'\circ\mathrm{pr}_{U'}$.
In this case $\nb = \nb_{f'}\boxtimes d$
and our discussion of the irregular Hodge filtration
also reduces to the product situation.
For simplicity,
we will assume that $\dim\Delta(f) =n$ in the rest of this section.
The general case then can be deduced easily.

\subsection{The Newton polyhedron filtration}

We define the \textit{Newton polyhegron filtration $\NPF^\lam(\nb)$}
of $(\Om^\bullet(*S),\nb)$ on $\torX$
similar to the filtration $F^\lam(\nb)$ for a good compactification $X$.
Again let $P$ be the pole divisor of $f$ on $\torX$.
Let
\begin{equation}\label{Eq:Def-NPF}
\NPF^\lam(\nb) :=
	\left[ \OO(\flr{-\lam P}) \xrightarrow{\nb} \OL^1(\flr{(1-\lam) P})
	\to\cdots\to \OL^p(\flr{(p-\lam) P}) \to\cdots \right]^{\geq\ceil{\lam}}.
\end{equation}
Notice that if the origin is contained in the interior of $\Delta(f)$,
then the morphism $f:U\to\AA$ is proper
and the filtration $\NPF^\lam(\nb)$ is indeed exhaustive.

To compute the hypercohomology of $\NPF^\lam(\nb)$,
first notice that on the toric variety $\torX$
the locally free sheaf $\OL^p$ is trivial for any $p$.
Indeed as an $\OO$-module it is globally generated by
\[ \bigwedge^p \left\{
	\frac{dx_1}{x_1},\cdots,\frac{dx_n}{x_n} \right\}. \]
On the other hand for $p\geq\lam$,
we have
\[ H^i\big(\torX,\OO(\flr{(p-\lam)P}\big) = 0
\quad\text{if $i\neq 0$} \]
and
$H^0\big(\torX,\OO(\flr{(p-\lam)P}\big)$ equals
the $\CC$-vector space generated by $\{x^\alpha\}$
where $\alpha$ runs over the lattice points
inside the dilated polyhedron $(p-\lam)\cdot\Delta(f)$
(see \cite[Prop.p.68 and Cor.p.74]{F}).
Thus one obtains
\begin{equation}\label{Eq:Toric-hyperH-H}
\HH^i\left(\torX,\NPF^\lam(\nb)\right) =
	H^i\left(\Gamma\left(\torX,\NPF^\lam(\nb)\right)\right)
\end{equation}
and this cohomology does not depend on the choice of $\torX$.

\begin{Thm}[Adophson-Sperber]\label{Ass:AS-Nagoya}\footnote{
The relations between the notations here and in \cite{AS-Nagoya}
are that
\[ \Gamma\left(\torX,\NPF^\lam(\nb)^l\right) = F_{-\lam}\hat{K}^l
\quad\text{and}\quad
\GrNP^\bullet(\nb) = (\bar{K}^\cdot,\bar{\delta}_{f,\alpha}). \]
The proof of (i)
is established in \cite[pp.70-73]{AS-Nagoya}
where the authors show
the quasi-isomorphism between the two complexes
$(\hat{K}^\cdot,\delta_{f,\alpha})$ and $(K_0^\cdot,\delta_{f,\alpha})$,
which correspond to
$\NPF^0(\nb)$ and $(\Om^\bullet(*S),\nb)$,
respectively.
For (ii),
see \cite[(4.3)]{AS-Nagoya},
cf.~\cite[Thm.2.14]{AS-Annals} and \cite[Th.2.8]{K}.
For (iii),
see {\cite[Thm.1.4 and Thm.4.1]{AS-Nagoya}}.
Notice that there is a typo in \cite[p.68]{AS-Nagoya}.
In line 18, the weight $(k/e)-l$ should be $(k/e)+l$.}
Suppose $\Delta(f) = n$ and
$f$ is non-degenerate with respect to $\Delta(f)$.
With notations as above, we have the following.
\begin{enumerate}
\item
For $\lam \leq0$
the inclusion $\NPF^\lam(\nb) \to (\Om^\bullet(*S),\nb)$
on $X$ is a quasi-isomorphism.
\item
$H^i\left(\Gamma\left(\torX,\GrNP^\lam(\nb)\right)\right) \neq 0$
only if $i = n$.
\item
Let $\mathrm{Vol}(f)$ be the usual Euclidean volume
of $\Delta(f)$ in $\RR^n$.
Then
\[ \dim \HdR^i(U,\nb) = \left\{\begin{array}{ll}
n!\cdot\mathrm{Vol}(f) & \text{if $i=n$} \\
0 & \text{otherwise}. \end{array}\right. \]
\end{enumerate}
\end{Thm}

Combined with \eqref{Eq:Toric-hyperH-H},
the above theorem implies the following.

\begin{Cor}\label{Cor:NPF-deRham-deg}
Suppose $\Delta(f) = n$ and
$f$ is non-degenerate with respect to $\Delta(f)$.
The spectral sequence attached to the filtration
$\NPF^\lam(\nb)$ on $\torX$
converges to $\HdR(U,\nb)$
and degenerates at the initial stage.
\end{Cor}

\subsection{The comparison}
As already mentioned,
the rational function $f$ on $\torX$ is not yet a morphism to $\PP$ in general
and hence $(\torX,S)$ is not a good compactification of $(U,f)$
for defining the irregular Hodge filtration.
This is because the zero divisor $Z$ and the pole divisor $P$ of $f$
intersect and one needs to perform blowups,
say $\pi: X \to \torX$,
in order to eliminate the indeterminacy.
However when $f$ is non-degenerate with respect to $\Delta(f)$,
we can say more.

\begin{Prop}\label{Prop:Geom-non-degen}
Suppose that $f$ is non-degenerate with respect to $\Delta(f)$.
Then on $\torX$ the zero divisor $Z$
and the support of the pole divisor $\red{P}$ of $f$
intersect transversally
and the intersections of $Z$ with various toric strata of $\red{P}$
are smooth.
\end{Prop}

\pf
A codimension $r$ toric stratum $D$ of $\red{P}$ is a dense torus
sitting in an irreducible component of the intersection
of certain irreducible components $D_1,\cdots,D_r$ of $S$.
Each $D_i$ corresponds to a ray in $N_\RR$,
which then corresponds to a face $\delta_i$ of $\Delta(f)$
(containing the exponents $\alpha\in\Delta(f)$
with most negative product with the direction of the ray).
A face $\delta$ in the intersection of $\delta_i$
then corresponds to $D$
and $f_\delta$ is the most singular term
of the function $f$ restricted to $D$
since those monomials in $f_\delta$
are among the terms in $f$ with the highest pole order along $D$.
We have $0\not\in\delta$
since otherwise $f$ has no pole along $D$.
Also the indeterminacy locus $Z\cap D$ on $D$
is exactly the zero set defined by $f_\delta =0$
(with variables along the $\delta$-direction).
Now the condition of emptiness of the solution of \eqref{Eq:Non-deg-cond}
(which becomes the usual Jocabian criterion after a change of variables)
exactly says that $Z\cap D$ is smooth,
which is what we want.
\qed
\bigskip

From now on
we assume that $f$ is non-degenerate with respect to $\Delta(f)$.
\medskip

We construct one particular $\pi:X\to\torX$ to obtain $f:X\to\PP$
as follows.
One picks an irreducible component $D$ of $Z\cap\red{P}$
and then take the blowup $X'$ along $D$.
If the exceptional divisor $E$ contributes to the pole
of $f$ on $X'$,
we perform the blowup along $E\cap Z'$
where $Z'$ is the zero divisor of $f$ on $X'$.
Continue this procedure until $f$ extends to a morphism to $\PP$
along the exceptional locus on $X^{(k)}$.
Let $Z^{(k)}$ and $P^{(k)}$ be the zero and pole divisors of $f$ on $X^{(k)}$.
Then one picks one irreducible component of $Z^{(k)}\cap\red{P^{(k)}}$
and performs a sequence of blowups again as above.
Repeating the procedure,
one then obtains the commutative diagram
\begin{equation}\label{Eq:X-torX}
\xymatrix{
X \ar[r]\ar@/^2pc/[rrrrr]^\pi\ar@/_1pc/[drrrrr]_f &
\cdots \ar[r] & X_2 \ar[r]^\eps & X_1 \ar[r] & \cdots \ar[r] &
\torX \ar@{-->}[d]^f\\
&&&&& \PP }
\end{equation}
where each step is the blowup
along a smooth irreducible component of the intersection of
the zero and pole divisors of $f$.
\medskip

Now $(X,X\setminus U)$ defines a good compactification of $(U,f)$
and we have the filtration $F^\lam(\nb)$ on $X$.
We shall show that there is a natural quasi-isomorphism
between $\RD\pi_*F^\lam(\nb)$ and $\NPF^\lam(\nb)$ on $\torX$
for each $\lam$.
Consequently they define the same filtration on $\HdR(U,\nb)$
and furthermore the Hodge to de Rham spectral sequence degenerates
in this case.

For this and to simplify the notations,
we consider the filtrations,
called $F_1^\lam(\nb)$ and $F_2^\lam(\nb)$,
of the twisted de Rham complexes
on $X_1$ and $X_2$, respectively
where $\eps:X_2\to X_1$ appears in the above sequence of blowups.
The two filtrations $F_i^\lam(\nb)$
are defined exactly as in \eqref{Eq:Def-NPF}
(which does not require the variety is toric).
Now notice that
for $X_2=X$
the filtration $F_2^\lam(\nb)$ is $F^\lam(\nb)$
for the good compactification $X$
while for $X_1=\torX$
the filtration $F_1^\lam(\nb)$ is the Newton polyhedron filtration
$\NPF^\lam(\nb)$ on the toric $\torX$.

\newcommand{\Bl}{\operatorname{Bl}}
We look at the local situation
over a point of the center of blowup in $X_1$.
Prop.\ref{Prop:Geom-non-degen} ensures the following.
We can take $X_1 = \Disc^n$
with coordinates
\[ \{x,y_1,\cdots,y_k,t_1,\cdots,t_l,z,\tau_1,\cdots,\tau_m\} \]
and $U = (\Disc^\circ)^{1+k+l}\times\Disc^{1+m}$
with the boundary $S_1 = (xyt)$.
The regular function on $U$ is
\[ f = \frac{z}{x^ey^r} = \frac{z}{x^ey_1^{r_1}\cdots y_k^{r_k}}
	\quad(\text{for some $e,r>0$}). \]
The center $\Xi$ of blowup $\eps$ is
given by $x=0=z$.

The blowup $X_2 \subset X_1\times\PP$
is given by the equation
\[ xu = zv
	\quad([u,v]\in\PP). \]
We use the notations in the illustration of
$X_2\xrightarrow{\eps}X_1$ below.

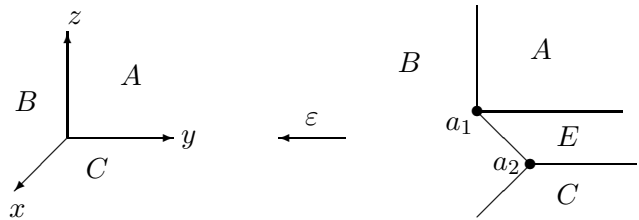
\begin{picture}(250,100)(-75,0)
\put(25,35){\vector(1,0){40}}
\put(25,35){\vector(0,1){40}}
\put(25,35){\vector(-1,-1){20}}
\put(3,5){$x$}
\put(25,78){$z$}
\put(68,33){$y$}
\put(45,55){$A$}
\put(5,45){$B$}
\put(32,20){$C$}

\put(180,45){\line(1,0){55}}
\put(180,45){\line(0,1){40}}
\put(180,45){\line(1,-1){20}}
\put(200,25){\line(-1,-1){20}}
\put(200,25){\line(1,0){40}}
\put(200,65){$A$}
\put(150,60){$B$}
\put(210,10){$C$}
\put(210,32){$E$}
\put(180,45){\circle*{4}}
\put(168,38){$a_1$}
\put(200,25){\circle*{4}}
\put(186,23){$a_2$}

\put(130,35){\vector(-1,0){25}}
\put(115,40){$\eps$}
\end{picture}

\noin Here the exceptional divisor $E$ is a split $\PP^1$-bundle
over the $y$-$t$-$\tau$-coordinate plane $\Xi$ of $X_1$.
Let
\[ e':=e-1. \]
The pole divisors of $f$ on $X_1$ and $X_2$
are given by
\[ P_1 = eA+rB
	\quad\text{and}\quad
P_2 = eA+e'E+rB, \]
respectively.
We have the information at the two points $a_1$ and $a_2$
in the table below
with $\bar{v} = \frac{v}{u}$ and $\bar{u} = \frac{u}{v}$.
\[\begin{array}{c|cc}
& a_1 & a_2 \\
\hline
\text{coordinates} & \{ \bar{v},z,y,t,\tau \} & \{ \bar{u},x,y,t,\tau \} \\
f & \frac{1}{\bar{v}^ez^{e'}y^r} & \frac{\bar{u}}{x^{e'}y^r}
\end{array}\]

Now with the index $\lam$ fixed,
we consider a sequence of new complexes
as follows.
Write $p = \ceil{\lam}$.
For $q = 0,1,\cdots,(n-p)$,
let $R_\lam(q) = (R_\lam(q)^\bullet,\nb)$
be the complex on $X_2$ given by
\begin{equation}\label{Eq:Toric-refine-def}
R_\lam(q)^{p+j} = \left\{\begin{array}{ll}
0 & \text{if $j< 0$} \\
\OL^{p+j}\big((\mu+je)A+(\mu+je)E+(\nu+jr)B\big) &
	\text{if $0\leq j\leq q$} \\
\OL^{p+j}\big((\mu+je)A+(\mu+je'+q)E+(\nu+jr)B\big) &
	\text{if $j>q$} \end{array}\right.
\end{equation}
where $\mu = \flr{(p-\lam)e}$ and $\nu = \flr{(p-\lam)r}$.
We have
\begin{eqnarray}\label{Eq:Trunc-R}
R_\lam(q-1)^{\leq p+q-1} &=& R_\lam(q)^{\leq p+q-1} \\
\label{Eq:Toric-pi-R}
\pi^*\left(F_1^\lam(\nb)\right) &\subset& R_\lam(n-p) \\
\label{Eq:Toric-refine}
R_\lam(-1) &:=& F_2^\lam(\nb) \subset R_\lam(0) \subset R_\lam(1)
	\subset\cdots\subset R_\lam(n-p).
\end{eqnarray}
Notice that $R_\lam(-1) = R_\lam(0)$
if $\flr{(p-\lam)e} = \flr{(p-\lam)e'}$.

\begin{Prop}\label{Prop:Toric-F-F}
With notations as above,
we have the following.
\begin{enumerate}
\item The quotient of the inclusion \eqref{Eq:Toric-pi-R}
is a complex with each component
equal to a direct sum
of the relative degree $(-1)$-invertible sheaves $\OO_{E/\Xi}(-1)$.
\item For $q=0,1,\cdots,n-p$,
the quotient $R_\lam(q)/R_\lam(q-1)$ is quasi-isomorphic to a direct sum
of $\OO_{E/\Xi}(-1)$ concentrated at degree $p+q$.
\end{enumerate}
\end{Prop}

\pf
For (i),
we have
\[ R_\lam(n-p)^j/ \pi^* F_1^\lam(\nb)^j \cong
	\OO_{E/\Xi}(-1)^{\binom{n-1}{j-1}} \]
for $j\geq p$ by a direct computation.

To understand the successive quotients of \eqref{Eq:Toric-refine},
we introduce one more complex.
Let $S_2=X_2\setminus U$.
For three integers $\rho,\eta,\xi$,
we let $(K_{\rho,\eta,\xi}^\bullet,\nb)$
be the subcomplex of $(\Om^\bullet(*S_2),\nb)$
on $X_2$
whose degree-$j$ term is given by
\[ K_{\rho,\eta,\xi}^j =
	\OL^j\left((\rho+je)A+(\eta+je')E+(\xi+jr)B\right). \]
One has
\begin{equation}\label{Eq:Rel-R-K}
R_\lam(q)^{\geq p+q} =
	\left(K_{\mu-pe,\mu-pe'+q,\nu-pr}^\bullet,\nb\right)^{\geq p+q}.
\end{equation}

\begin{Lemma}
The inclusion
\[ \left( K_{\rho,\eta,\xi}^\bullet,\nb \right) \subset
	\left( K_{\rho,\eta+1,\xi}^\bullet,\nb \right) \]
is a quasi-isomorphism of complexes on $X_2$
for any $\rho,\eta,\xi\in\ZZ$
with the condition that $\eta\geq 0$ if $e'=0$.
\end{Lemma}

\pf
At the point $a_1$,
one only needs to consider the case where $e=1$,
thanks to Prop.\ref{Prop:Change-Support}.
In this case,
we have the exterior product decomposion
\[ (K_{\rho,\eta,\xi}^\bullet,\nb) = (K_{\rho,\xi}^\bullet,\nb') \boxtimes
\left[ \frac{1}{z^\eta}\OO_z \xrightarrow{d} \frac{1}{z^\eta}\OL_z^1 \right] \]
where $K_{\rho,\xi}^\bullet$ is defined similarly
as the definition of $K_{\rho,\eta,\xi}^\bullet$ above
but now on the coordinates $\{\bar{v},y,t,\tau\}$
for the connection $\nb'$ attached to $1/(\bar{v}y^r)$.
The assertion then follows from the fact that
\[ \left[ \OO_z \xrightarrow{d} \OL_z^1 \right] \to
	\left[ \frac{1}{z^i}\OO_z \xrightarrow{d} \frac{1}{z^i}\OL_z^1 \right] \]
is a quasi-isomorphism for any $i\geq 0$.

The case for points between $a_1$ and $a_2$
is similar.

At the point $a_2$,
let $\OO$ be the coordinate ring
and $\bar\OO = \OO/x\OO$.
One has to check the exactness of
\begin{equation}\label{Eq:a2-j-exact}\begin{split}
\frac{\OL^{j-1}((\eta+1-e')E+(\xi-r)B)}{\OL^{j-1}((\eta-e')E+(\xi-r)B)}
&\xrightarrow{\nb_{j-1}}
	\frac{\OL^j((\eta+1)E+\xi B)}{\OL^j(\eta E+\xi B)} \\
&\xrightarrow{\nb_j}
	\frac{\OL^{j+1}((\eta+1+e')E+(\xi+r)B)}{\OL^{j+1}((\eta+e')E+(\xi+r)B)}
\end{split}\end{equation}
where now the connection is the $\bar\OO$-linear map
given by the left cup product with
\[ \frac{1}{x^{e'}y^r}\left(d\bar{u}-\bar{u}\frac{e'dx}{x}-
	\bar{u}\sum_{i=1}^k\frac{r_idy_i}{y_i}\right)
= \frac{\bar{u}}{x^{e'}y^r}\left(\frac{d\bar{u}}{\bar{u}}-\frac{e'dx}{x}-
	\sum_{i=1}^k\frac{r_idy_i}{y_i}\right). \]
First suppose that $e'\geq 1$.
Let
\[ \Lambda_i := \bar\OO\cdot\bigwedge^i
	\left\{d\bar{u},\frac{dx}{x},\frac{dy}{y},\frac{dt}{t},d\tau\right\} \]
be the $\bar\OO$-module
generated by $i$-forms.
Notice that the complex
\[ \frac{1}{\bar{u}^2x^{\eta+1-2e'}y^{\xi-2r}}\Lambda_{j-2}
\xrightarrow{\nb_{j-2}}
	\frac{1}{\bar{u}x^{\eta+1-e'}y^{\xi-r}}\Lambda_{j-1}
\xrightarrow{\nb_{j-1}}
	\frac{1}{x^{\eta+1} y^\xi}\Lambda_j \xrightarrow{\nb_j}
	\frac{\bar{u}}{x^{\eta+1+e'}y^{\xi+r}}\Lambda_{j+1}, \]
being isomorphic to the Koszul complex
associated with
\[ \{\gamma_\infty,\gamma_0,\gamma_i\}_{i=1}^k
\quad\text{corresponding to}\,
\frac{\bar{u}}{x^{e'}y^r}\left\{\frac{d\bar{u}}{\bar{u}},\frac{-e'dx}{x},
	\frac{-r_idy_i}{y_i}\right\}_{i=1}^k, \]
is exact.
Thus if $\nb_j(\alpha)=0$
for some $\alpha$ in our degree-$j$ piece,
there exists a
\[ \beta = d\bar{u}\wedge\beta_1+\beta_2
	\quad\text{with}\, \left\{\begin{array}{l}
\beta_1 \in \frac{1}{\bar{u}^2x^{\eta+1-e'}y^{\xi-r}}\bar\OO\cdot
	\bigwedge^{j-2}\left\{\frac{dx}{x},\frac{dy}{y},\frac{dt}{t},d\tau\right\} \\
\beta_2 \in \frac{1}{\bar{u}x^{\eta+1-e'}y^{\xi-r}}\bar\OO\cdot
	\bigwedge^{j-1}\left\{\frac{dx}{x},\frac{dy}{y}\frac{dt}{t},d\tau\right\}
\end{array}\right. \]
such that $\nb_{j-1}(\beta)=\alpha$.
By subtracting $\nb_{j-2}(x^{e'}y^r\beta_1)$ to $\beta$,
we may assume that $\beta_1=0$.
Then the part
$x^{-e'}y^{-r}d\bar{u}\wedge\beta_2$ of $\alpha=\nb_{j-1}(\beta)$
does not have $\bar{u}$ in the denominator,
and hence neither does $\beta_2$.
Therefore
$\beta\in\frac{1}{x^{\eta+1-e'}y^{\xi-r}}\Lambda_{j-1}$
and \eqref{Eq:a2-j-exact} is exact.

The case $e'=0$ is similar.
\qed
\bigskip

Now back to the proof of (ii) in Prop.\ref{Prop:Toric-F-F}.
Let $\mu = \flr{(p-\lam)e}$ and $\nu = \flr{(p-\lam)r}$.
Let $a=\mu-pe$, $b=\mu-pe'+q$ and $c=\nu-pr$.
By the relations \eqref{Eq:Trunc-R} and \eqref{Eq:Rel-R-K}
and the previous lemma,
we have
\begin{equation}\label{Eq:Q-K-I}
\frac{R_\lam(q)}{R_\lam(q-1)} \xleftarrow{\sim}
\ker\left\{\frac{K_{a,b,c}^{p+q}}{K_{a,b-1,c}^{p+q}}\xrightarrow{\nb}
	\frac{K_{a,b,c}^{p+q+1}}{K_{a,b-1,c}^{p+q+1}}\right\}
\xleftarrow{\sim}
\Image\left\{\frac{K_{a,b,c}^{p+q-1}}{K_{a,b-1,c}^{p+q-1}}\xrightarrow{\nb}
\frac{K_{a,b,c}^{p+q}}{K_{a,b-1,c}^{p+q}}\right\}.
\end{equation}
(Notice that if $e'=0$ and $b=0$,
then we have $q=\mu=0$
and $R_\lam(-1) = R_\lam(0)$.
So there is nothing to prove.)

Let $\eta = \mu+qe$ and $\xi = \nu+qr$.
Away from $u=0$ (neighborhood of $a_1$),
the $\OO_E$-module $K_{a,b,c}^{p+q-1}/K_{a,b-1,c}^{p+q-1}$
is generated by
$\om_1$ of the four types listed in the table
\[\begin{array}{c|c}
\bar{v}^{\eta-e}z^{\eta-e'}y^{\xi-r}\cdot\om_1 &
x^{\eta-e'}y^{\xi-r}\cdot\om_2 \\
\hline
\zeta_1\cdots\zeta_{p+q-1} &
\zeta_1\cdots\zeta_{p+q-1} \\
\frac{d\bar{v}}{\bar{v}}\zeta_1\cdots\zeta_{p+q-2} &
-\left(e'\frac{dx}{x}+\sum_{i=1}^k r_i\frac{dy_i}{y_i}\right)
	\zeta_1\cdots\zeta_{p+q-2} \\
\frac{dz}{z}\zeta_1\cdots\zeta_{p+q-2} &
\left( e\frac{dx}{x}+\sum_{i=1}^k r_i\frac{dy_i}{y_i}\right)
	\zeta_1\cdots\zeta_{p+q-2} \\
\frac{d\bar{v}}{\bar{v}}\frac{dz}{z}
	\zeta_1\cdots\zeta_{p+q-3} &
\left( \frac{dx}{x}\sum_{i=1}^k r_i\frac{dy_i}{y_i}\right)
	\zeta_1\cdots\zeta_{p+q-3}
\end{array}\]
where
\[ \{\zeta_j\} = \left\{\frac{dy}{y},\frac{dt}{t},d\tau\right\}. \]
One checks that we have
\[ \bar{u}\cdot\nb(\om_1) \equiv \nb(\om_2)
	\pmod{K_{a,b-1,c}^{p+q} = \OL^{p+q}\big(\eta A+(\eta-1)E+\xi B\big)} \]
where $\om_2$ are the corresponding forms lying away from $v=0$
(neighborhood of $a_2$) listed above.
This equation shows that the last term of \eqref{Eq:Q-K-I}
is a direct sum of $\OO_{E/\Xi}(-1)$.
\qed

\begin{Thm}
Consider the pair $(U,f)$
where $U$ is a torus of dimension $n = \dim\Delta(f)$
and $f$ is non-degenerate with respect to $\Delta(f)$.
Then
the irregular Hodge filtration coincides
with the filtration induced by $\NPF^\lam(\nb)$
on any smooth toric compactification $\torX$
with simple normal crossing boundary $\torX\setminus U$,
and the irregular Hodge to de Rham spectral sequence degenerates
at the initial stage.
\end{Thm}

\pf
We choose a good compactification $X$
with $\pi:X\to\torX$
as constructed in \eqref{Eq:X-torX}.
By Prop.\ref{Prop:Toric-F-F},
we have a natural quasi-isomorphism
between $\RD\pi_*F^\lam(\nb)$ and $\NPF^\lam(\nb)$
on $\torX$ for any $\lam$.
The assertions now follow from
Thm.\ref{Ass:AS-Nagoya} and Cor.\ref{Cor:NPF-deRham-deg}.
\qed

\appendix
\newcommand{\FD}{\mathfrak{F}}
\section{Comparison with Deligne's definition}

In this appendix
we recall Deligne's definition of the irregular Hodge filtration
in the curve case in \cite{DMR}
and show that it induces the same filtration as ours
in the de Rham cohomology.
\medskip

Consider the pair $(U,f)$ where $U$ is a smooth curve.
Let $X$ be the smooth completion of $U$
with boundary $S:= X\setminus U$.
Let $\nb = d+df$ and write $P=$ the pole divisor of $f$ on $X$
as before.

Deligne then defines inductively
an exhaustive and separated filtration $\FD^\lam$
of the two-term complex $(\Om_X^\bullet(*S),\nb)$ by letting
\[ \FD^\lam (\Om_X^\bullet(*S),\nb) =
	\left[\FD^\lam\OO_X(*S) \xrightarrow{\nb}
		\FD^\lam\Om_X^1(*S)\right] \]
where
\begin{eqnarray*}
\FD^\lam \OO_X(*S) &=& \left\{\begin{array}{ll}
	0 & \text{if $\lam > 0$} \\
	\OO_X(S-\ceil{\lam P}) & \text{if $-1 < \lam \leq 0$} \\
	\left(\FD^{\lam+1}\OO_X(*S)\right)(S+P) &
		\text{if $\lam \leq -1$}\end{array}\right.
	\\ \\
\FD^\lam \Om_X^1(*S) &=&
	\Om_X^1\tensor{\OO_X}\left(\FD^{\lam-1}\OO_X(*S)\right).
\end{eqnarray*}

Define a subcomplex
$\Om_X^\bullet(\log_\FD\nb)$ of $(\Om_X^\bullet(*S),\nb)$
to be the two-term complex
\[\left[\xymatrix{
\ker\left\{ \FD^0\OO_X(*S) \xrightarrow{\nb} \Gr_\FD^0\Om^1_X(*S)\right\}
	\ar[r]^(.68)\nb & \FD^{0+}\Om^1_X(*S) }\right] \]
equipped with the induced filtration $\FD^\lam$.
Then $\FD^\lam\Om_X^\bullet(\log_\FD\nb)$ is non-trivial
only if $0\leq\lam\leq 1$.
We call $\Om_X^\bullet(\log_\FD\nb)$
the logarithmic subcomplex of $(\Om_X^\bullet(*S),\nb)$;
it is a complex of coherent sheaves on $X$,
filtered by coherent subcomplexes.
The context of the irregular Hodge theory over curves
is summarized as the following.

\begin{Thm}[Deligne]
With notations as above,
we have the following.
\begin{enumerate}
\item The natural inclusion
$(\Om_X^\bullet(\log_\FD\nb),\nb,\FD) \to (\Om_X^\bullet(*S),\nb,\FD)$
is a quasi-isomorphism of filtered complexes on $X$.
\item For each $\lam$,
the map $\HH(X,\FD^\lam) \to \HdR(U,\nb)$
induced by the inclusion of complexes is injective
(i.e.~the spectral sequence associated with the filtration $\FD$
degenerates at the initial $E_1$ stage).
\end{enumerate}
\end{Thm}

We remark again
that the construction can be generalized to
exponential twists of unitary regular connections of any ranks
over the curve $U$
and the corresponding statements as above continue to hold
in the general case.
\medskip

Now let us compare the two filtrations $F^\bullet(\nb)$ and $\FD^\bullet$.
First we clearly have $F^\lam(\nb) \subset \FD^\lam$
for any $\lam \in \RR$.
On the other hand, one readily observes that
the two corresponding logarithmic filtered complexes
$\Om_X^\bullet(\log_F\nb)$ and $\Om_X^\bullet(\log_\FD\nb)$
are exactly the same subcomplex of $(\Om_X^\bullet(*S),\nb)$.
Thus we obtain the following statement.

\begin{Prop}
In the curve case,
the two filtrations $F^\bullet(\nb)$ and $\FD^\bullet$
induce the same filtration
on the twisted de Rham cohomology $\HdR(U,\nb)$.
\end{Prop}

\end{document}